\documentclass[11pt]{article}
\usepackage{amsmath,amsthm,amssymb,color,epic,eepic,
cite}
\newlength{\minitwocolumn}
\newcommand{\Z}{{\Bbb Z}} 
\newcommand{\C}{{\Bbb C}} 
\newcommand{\FF}{{\Bbb F}} 
\newcommand{\K}{{\Bbb K}} 

\newcommand{\cD}{{\cal D}}

\newcommand{\cR}{{\cal R}}

\newcommand{\cQ}{{\cal Q}}

\newcommand{\hL}{\widehat{L}}

\newcommand{\la}{\lambda}

\newcommand{\al}{\alpha}

\newcommand{\vep}{\varepsilon}

\newcommand{\bH}{\bar{H}}

\newcommand{\ha}{{\alpha}}
\newcommand{\hb}{{\beta}}
\newcommand{\hf}{\widehat{f}}

\newcommand{\hV}{\widehat{V}}
\newcommand{\hW}{\widehat{W}}
\newcommand{\hOmega}{\widehat{\Omega}}

\newcommand{\noi}{{\noindent}}
\newcommand{\nn}{{\nonumber}}
\newcommand{\bea}{\begin{eqnarray}}
\newcommand{\ena}{\end{eqnarray}}
\newcommand{\beit}{\begin{itemize}}
\newcommand{\enit}{\end{itemize}}

\newcommand{\be}{\begin{eqnarray*}}
\newcommand{\en}{\end{eqnarray*}}
\newcommand{\lb}[1]{\label{#1}}

\newcommand{\ds}[1]{{\displaystyle #1 }}

\newcommand{\End}{{\rm End}}
\newcommand{\rank}{{\rm rank}}

\newcommand{\id}{{\rm id}}

\def\infq4p#1{{(#1;q^4,p)_\infty}}

\newcommand{\hLp}{\widehat{L}^+}

\newcommand{\tot}{\widetilde{\otimes}}

\newcommand{\mmatrix}[1]{\begin{matrix} #1 \end{matrix}}
\newcommand{\mat}[1]{\left(\mmatrix{#1}\right)}

\font\teneufm=eufm10
\font\seveneufm=eufm7
\font\fiveeufm=eufm5
\newfam\eufmfam
\textfont\eufmfam=\teneufm
\scriptfont\eufmfam=\seveneufm
\scriptscriptfont\eufmfam=\fiveeufm

\let\goth\frak
\newcommand{\slth}{\widehat{\goth{sl}}_2}
\newcommand{\slt}{\goth{sl}_2}
\newcommand{\gsl}{\goth{sl}}
\newcommand{\slnh}{\widehat{\goth{sl}}_N}
\newcommand{\sln}{\goth{sl}_N}
\newcommand{\g}{\goth{g}}

\newcommand{\uq}{U_q\bigl(\slth\bigr)}
\newcommand{\Aqp}{{\cal A}_{q,p}}

\newcommand{\Bqla}{{{\cal B}_{q,\lambda}}}
\newcommand{\Uqp}{U_{q,p}}

\newcommand{\h}{H}
\newcommand{\hh}{\goth{h}}

\font\fourteeneufm=eufm10 scaled\magstep2    
\font\seventeeneufm=eufm10 scaled\magstep3   
   
\newcommand{\slthbig}{\widehat{\mbox{\fourteeneufm sl}}_2}  
\newcommand{\slthBig}{\widehat{\mbox{\seventeeneufm sl}}_2} 
\makeatletter
\@addtoreset{equation}{section}
\makeatother

\newtheorem{thm}{Theorem}[section]
\newtheorem{prop}[thm]{Proposition}
\newtheorem{lem}[thm]{Lemma}
\newtheorem{cor}[thm]{Corollary}

\newtheorem{dfn}[thm]{Definition}

\begin{document}
\bibliographystyle{unsrt}

\

\vspace{-1cm}
\begin{center}
{\bf\Large  Elliptic Quantum Group $U_{q,p}(\slthBig)$, 
 Hopf Algebroid Structure 
 and
  Elliptic Hypergeometric Series\\[7mm] }
{\large  Hitoshi Konno}\\[6mm]
{\it Department of Mathematics, Graduate School of Science, 
\\Hiroshima University, Higashi-Hiroshima 739-8521, Japan\\
       konno@mis.hiroshima-u.ac.jp}\\[7mm]
\end{center}

\begin{abstract}
\noindent 
We propose a new realization of the 
elliptic quantum group equipped with the $\h$-Hopf algebroid structure 
on the basis of the elliptic algebra $U_{q,p}(\slth)$.  
The algebra $U_{q,p}(\slth)$ has a constructive definition   
in terms of the Drinfeld generators of the 
quantum affine algebra $U_q(\slth)$ 
and a Heisenberg algebra. 
This yields a systematic construction of  both finite and infinite-dimensional 
dynamical representations and  their parallel structures to $U_q(\slth)$.  
In particular we give a classification theorem of
 the finite-dimensional irreducible 
pseudo-highest weight representations stated  
in terms of an elliptic analogue of the 
Drinfeld polynomials. We also investigate a structure of 
the tensor product of two evaluation representations 
and derive an elliptic analogue of the Clebsch-Gordan 
coefficients. We show that it  
 is expressed by using the very-well-poised balanced elliptic 
hypergeometric series ${}_{12}V_{11}$. 
\end{abstract}
\nopagebreak

\section{Introduction}
In this paper we revisit the elliptic algebra  $U_{q,p}(\slth)$ and 
study its coalgebra structure. 
The algebra $U_{q,p}(\slth)$ was introduced in \cite{Konno}
 as an elliptic analogue of 
 the quantum affine algebra  $U_q(\slth)$ in the Drinfeld realization\cite{Drinfeld}.  
It was realized in  \cite{JKOS2} that 
$U_{q,p}(\slth)$ is constructively defined
 by using the Drinfeld 
 generators of $U_q(\slth)$ and 
the Heisenberg algebra  $\{P, e^Q\}$. 
This construction was generalized to   
the elliptic algebra $U_{q,p}(\g)$ of all types of untwisted affine 
Lie algebras $\g$ \cite{JKOS2} and of 
the twisted type $A^{(2)}_2$ \cite{KK04}.  

It was also realized that 
$U_{q,p}(\g)$ has an interesting relation to the deformed 
coset Virasoro/$W$ algebras\cite{SKAO,FF,AKOS,FR}.  
Namely, the level one ($c=1$) elliptic currents of $U_{q,p}(\g)$ 
are identified with the screening currents of 
the deformed $W$ algebras for $\g=\slnh$ \cite{Konno,JKOS2,KK03} 
and for 
$\g=A^{(2)}_2$ \cite{KK04}. This observation led us to a conjecture that the 
elliptic currents $E_i(u)$ and $F_i(u)$ of $U_{q,p}(\g)$ $(i=1,2,\cdots,{\rm rank}\bar{\g})$ 
define the deformation 
of the Virasoro/$W$ algebra associated with 
the coset $\g\oplus \g/\g$ \cite{Konno,JKOS2}. 

A study of coalgebra structure on $U_{q,p}(\g)$ was far from  
straightforward. 
We constructed the $L$ operators 
in terms of the elliptic currents and derived  
the $RLL$ relation 
for  the cases $\g=\slnh$  \cite{JKOS2,KK03} and  
$\g=A^{(2)}_2$ \cite{KK04}. 
However it turned out that a naive 
Fadeev-Reshetikhin-Sklyanin-Takhtajan (FRST) construction
\cite{FRT,Sklyanin} does not work 
due to the dynamical shift appearing in the $R$ matrices.  
Instead we obtained a connection to the quasi-Hopf algebra  
$\Bqla(\g)$ \cite{JKOS2}. That is 
the isomorphism $U_{q,p}(\g)\cong \Bqla(\g)\otimes \{P_i,e^{Q_i}\}$ 
$(i=1,2,\cdots,\rank \bar{\g})$ as an associative algebra, 
where $\{P_i,e^{Q_i}\}$ denotes a Heisenberg algebra.

The quasi-Hopf algebra $\Bqla(\g)$ (the face type)  
was introduced by Jimbo, Konno, Odake and Shiraishi \cite{JKOS} 
motivated by the works of Drinfeld\cite{DrinfeldQH}, Babelon, Bernard 
and Billey \cite{BBB} and Fr{\o}nsdal \cite{Fronsdal}. 
At the same time, we introduced the vertex type quasi-Hopf
 algebra  $\Aqp(\slnh)$. 
Both $\Aqp(\slnh)$ and $\Bqla(\g)$ are isomorphic to the corresponding 
quantum affine algebras 
$U_q(\g)$ as associative algebras, but their coalgebra structures 
are deformed from $U_q(\g)$ 
 by the twistors $E(r)$ and $F(\la)$, respectively. 
By twisting the objects in $U_q(\g)$, such as the comultiplication, the 
universal $R$ matrices and the vertex operators, 
we can derive their quasi-Hopf algebra counterparts\cite{JKOS}. 
Then the relation $U_{q,p}(\g)\cong \Bqla(\g)\otimes \{P_i,e^{Q_i}\}$ 
 allows us to derive  the $U_{q,p}(\g)$ counterparts from those of 
the quasi-Hopf algebra $\Bqla(\g)$. Such a strategy led us to an 
extension of the algebraic analysis scheme of trigonometric solvable lattice models 
\'{a} la Jimbo and Miwa \cite{JM} to  the face type elliptic
 models\cite{Konno,JKOS2,KK03,KK04,KKW}. 

However the tensor product with the Heisenberg algebras breaks down 
the quasi-Hopf algebra structure, so that $U_{q,p}(\g)$ is not a 
quasi-Hopf algebra. 
Moreover, the quasi-Hopf algebra itself has a disadvantage that 
its coalgebra structure is not suitable for a practical calculation 
due to a complication arising from the same twist procedure.  
This is a serious defect, for example, to develop the 
representation theory of the elliptic quantum groups  
and their  harmonic analysis. 

The aim of this paper is to  
show that a relevant coalgebra 
structure of $U_{q,p}(\g)$ is an $\h$-Hopf algebroid 
 and to formulate a new elliptic quantum group, which complements the 
disadvantage of the quasi-Hopf algebra. 
In this paper we consider the case $\g=\slth$. 
The cases of other affine Lie algebra types  will be 
discussed in future publications. 

The $\h$-Hopf algebroid was introduced by Etingof and Varchenko\cite{EV1,EV2} 
motivated by the works of Felder and Varchenko\cite{Felder,FV }. 
Some additional structures were given by Koelink and 
Rosengren\cite{KR,Rosengren}. A good review of this subject can be found in \cite{vN}.  
In \cite{EV1,TV,KNR}, it was applied to a formulation of Felder's elliptic quantum group 
$E_{\tau,\eta}(\slt)$ by using the generalized  
FRST construction on the basis of the 
$RLL$ relation associated with the elliptic dynamical 
$R$ matrix. 
 Another formulation of $E_{\tau,\eta}(\slt)$ 
as a quasi-Hopf algebra was studied 
by Enriquez and Felder\cite{EF}. 
 A similar Hopf algebroid structure was 
introduced by Lu \cite{Lu} and Xu \cite{Xu}. 
In \cite{Xu}, Xu also studied the 
algebra $\cD\otimes U_q(\bar{\g})$, 
where $\cD$ denotes the algebra of meromorphic 
differential operators on $\bar{\hh}^*$. His algebra is similar  
to $U_{q,p}(\g)$, but his $\bar{\g}$ is a finite-dimensional 
simple Lie algebra.  

Our formulation is based on  the fact  
 that the $RLL$ relation for $U_{q,p}(\slth)$ obtained in \cite{JKOS2} 
is identical to a central extension of the one for Felder's elliptic quantum group. 
This enables us to apply the generalized FRST construction to our case 
with a modification due to a central extension. 
The main  modification is that we use both 
the commuting subalgebra $H$ of $U_{q,p}(\slth)$ and 
the additive Abelian group $\bH^*\subset H^*$ appropriately in the formulation.
Here $H$ contains  the central element $c$, whereas $\bH$ does not.  
  We consider the field of meromorphic 
functions on $H^*$, but we use $\bH^*$ to define the bigrading structure of 
$U_{q,p}(\slth)$.   
As a result we obtain a new face type elliptic quantum group $U_{q,p}(\slth)$
 as an $\h$-Hopf algebroid, which 
  is realized in terms of 
 the Drinfeld generators and has the central extension. 
 We also show that the coalgebra structure is enough simple for 
a practical calculations. 

In comparison with the previous formulations\cite{FV,EV1,TV,KNR}, 
$U_{q,p}(\slth)$ has the advantage    
that it has a constructive definition in terms of the Drinfeld generators of 
 $U_q(\slth)$. This allows a systematic derivation of both finite and 
infinite-dimensional representations  of $U_{q,p}(\slth)$ from those of 
$U_q(\slth)$ and their parallel structures to $U_q(\slth)$. 
As an example, we study a classification theorem of the 
finite-dimensional irreducible pseudo-highest weight representations 
and make a statement in terms of an elliptic analogue of 
the Drinfeld polynomials. 
This gives an elliptic analogue of the works by Drinfeld\cite{Drinfeld} 
and by Chari and Pressley\cite{CP}. 
In addition, 
we investigate a submodule structure of the 
tensor product of two evaluation modules. 
We obtain the singular vectors explicitly and derive an 
elliptic analogue of the Clebsch-Gordan coefficients. 
We show that the coefficients are given by the terminating 
 very-well-poised balanced elliptic hypergeometric series ${}_{12}V_{11}$, which was 
introduced by Frenkel and Turaev \cite{FT} on the basis of the work by Date, Jimbo, Kuniba, 
Miwa and Okado\cite{DJKMO},   
and extensively studied by Spiridonov and 
Zhedanov \cite{S,SZ}. 
This provides the alternative 
to the representation theoretical derivation of  ${}_{12}V_{11}$
 by Koelink, van Norden and Rosengren \cite{KNR}. 
 In \cite{KNR}, ${}_{12}V_{11}$ 
was obtained as matrix elements of 
a co-representation of Felder's elliptic quantum group. 

In the separate paper \cite{Konno07}, we discuss  
 a free field representation of 
the infinite-dimensional highest weight representations of 
$U_{q,p}(\slth)$ and derive the vertex operators as intertwining operators of such 
$U_{q,p}(\slth)$-modules.  The resultant vertex operators coincide with 
those obtained indirectly in \cite{JKOS2} on the basis of the 
quasi-Hopf algebra structure of $\Bqla(\slth)$.    
This indicates a consistency of our 
  $\h$-Hopf algebroid structure on $U_{q,p}(\slth)$ 
  even in the case with non-zero central element.    
We hence establish the extension of the algebraic analysis scheme 
 to the fusion RSOS model on the basis of $U_{q,p}(\slth)$.

This paper is organized as follows. In Sect.2, we give a definition of the 
elliptic algebra $U_{q,p}(\slth)$ and review some results on the $RLL$ relation. 
In Sect.3, we recall the definition of $\h$-Hopf algebroid from \cite{EV1,EV2,KR}. 
Then we define an $\h$-Hopf algebroid structure on $U_{q,p}(\slth)$ and formulate it 
as an elliptic quantum group. In Sect.4, after a summary of basic facts on 
dynamical representations, we consider finite-dimensional representations of 
$U_{q,p}(\slth)$. In particular, we introduce an elliptic analogue of the Drinfeld polynomial 
and state a criterion for the finiteness of irreducible pseudo-highest weight 
representation of $U_q(\slth)$. We also investigate a 
submodule structure of the tensor product of two evaluation representations and derive 
 an elliptic analogue of the Clebsch-Gordan coefficients. 
Sect.5 is devoted to a discussion on the trigonometric, non-affine and non-dynamical 
limits of the results. 
In Appendix A, we give a list of the commutation relations of the $L$ operator elements. 
Appendix B is devoted to a proof of Theorem \ref{thmbetavs}.

\section{The Elliptic Algebra $U_{q,p}(\slthbig)$}
In this section, we give a definition of the elliptic algebra 
$\Uqp(\slth)$ in terms of the Drinfeld generators of the quantum affine algebra $U_q(\slth)$ and the Heisenberg algebra $\{P,e^Q\}$. 
We then recall some basic facts 
on $\Uqp(\slth)$ from \cite{Konno,JKOS2}.


\subsection{Quantum Affine Algebra $\K[U_q(\slth)]$}\lb{app:1.1}

Throughout this paper we fix a complex number $q$ such that $q\not=0, |q|<1$. 
\begin{dfn}\cite{Drinfeld}\lb{defUq}
For a field  $\K$, the quantum affine algebra $\K[U_q(\slth)]$ 
in the Drinfeld realization is an associative algebra over $\K$ 
generated by the standard Drinfeld generators 
$a_n\ (n\in \Z_{\not=0})$, $x_n^\pm\ (n\in \Z)$, 
$h$,  $c$ and $d$. 
The defining relations are given as follows.
\be
&&c :\hbox{ central },\nn\\
&& [h,d]=0,\quad [d,a_{n}]=n a_{n},\quad 
[d,x^{\pm}_{n}]=n x^{\pm}_{n}, \nn\\
&&[h,a_{n}]=0,\qquad [h, x^\pm(z)]=\pm 2 x^{\pm}(z),\nn\\
&&
[a_{n},a_{m}]=\frac{[2n]_{q}[c n]_{q}}{n}
q^{-c|n|}\delta_{n+m,0},\nn
\\&&
[a_{n},x^+(z)]=\frac{[2n]_{q}}{n}q^{-c|n|}z^n x^+(z),\nn\\
&&
[a_{n},x^-(z)]=-\frac{[2n]_{q}}{n} z^n x^-(z),
\lb{uqslth}
\en
\be
&&(z-q^{\pm 2}w)
x^\pm(z)x^\pm(w)= (q^{\pm 2}z-w) x^\pm(w)x^\pm(z),\nn
\\
&&[x^+(z),x^-(w)]=\frac{1}{q-q^{-1}}
\left(\delta\bigl(q^{-c}\frac{z}{w}\bigr)\psi(q^{c/2}w)
-\delta\bigl(q^{c}\frac{z}{w}\bigr)\varphi(q^{-c/2}w)
\right).\nn
\en
where we use $[n]_q=\frac{q^n-q^{-n}}{q-q^{-1}}$, $\delta(z)=\sum_{n\in\Z}z^n$ and the Drinfeld currents defined by 
\be
&&x^\pm(z)=\sum_{n\in \Z}x^\pm_{n} z^{-n},\\
&&\psi(q^{c/2}z)=q^{h}
\exp\left( (q-q^{-1}) \sum_{n>0} a_{n}z^{- n}\right),\quad 
\varphi(q^{-c/2}z)=q^{-h}
\exp\left(-(q-q^{-1})\sum_{n>0} a_{-n}z^{ n}\right).
\en

We also denote by $\K[U'_{q}(\slth)]$ the subalgebra of $\K[U_{q}(\slth)]$ 
generated by the same generators as $\K[U_{q}(\slth)]$ except for  $d$
 excluded.   
\end{dfn}

In the later sections, we use the symbols 
$\psi_n$ and $\phi_{-n}$ $(n\in \Z_{\geq 0})$ defined by
\be
&& \psi(q^{c/2}z)=\sum_{n\geq 0}\psi_n z^{-n},\qquad
 \varphi(q^{c/2}z)=\sum_{n\geq 0}\phi_{-n} z^{n}.
\en

 Let $\al_1$ and $\bar{\Lambda}_1=\frac{1}{2}\al_1$ be 
the simple root and  the fundamental weight of $\gsl(2,\C)$, respectively. 
We set $\bar{\hh}=\C h$, $\cQ=\Z \al_1$ and $\bar{\hh}^*=\C \bar{\Lambda}_1$. 
We denote by $<\ ,\ >$ the paring of $\bar{\hh}$ and 
$\bar{\hh}^*$ given by $<\bar{\Lambda}_1,h>=1$. 

\subsection{Definition of the Elliptic Algebra $U_{q,p}(\slth)$}\lb{app:1.3}

Let $r$ be a generic complex number. We set $r^*=r-c$, $p=q^{2r}$ and $p^*=q^{2r^*}$. 
We define the  Jacobi theta functions $[u]$ and $[u]^*$ by
\be
&&[u]=\frac{q^{\frac{u^2}{r}-u}}{(p;p)_\infty^3}\Theta_p(q^{2u}),\qquad [u]^*=[u]|_{r\to r^*},\\
&&\Theta_p(z)=(z;p)_\infty(p/z;p)_\infty(p;p)_\infty,
\en
where
\be
&&(z;p_1,p_2,\cdots, p_m)_\infty=\prod_{n_1,n_2,\cdots,n_m=0}^\infty(1-zp_1^{n_1}p_2^{n_2}\cdots p_m^{n_m}).
\en
Setting $p=e^{-\frac{2\pi i}{ \tau}}$, $[u]$ satisfies the quasi-periodicity $[u+r]=-[u]$, $[u+r\tau]=e^{-\pi i (2u/r+\tau)}[u]$.

Let $\{P,Q\}$ be a Heisenberg algebra commuting with $\C[U_q(\slth)]$ and 
satisfying 
\bea
&&[P, {Q}]=-1. \lb{Heisenberg}
\ena
We set $\h=\C P\oplus \C  r^*$ and  $\h^*=\C Q\oplus \C \frac{\partial}{\partial r^*}$. 
We denote by the same symbol $<\ ,\ >$ as the above the pairing of $\h$ and $\h^*$ 
defined by  
\be
&&<Q,P>=1=<\frac{\partial}{\partial r^*},r^*>,
\en
the others are zero. We regard $H^*\oplus H$ as a Heisenberg algebra by
\be
[x,y]=<x,y>.
\en

We also consider the Abelian group $\bar{\h}^{*}= \Z Q$. We have the isomorphism  
$\phi:\cQ\to \bar{\h}^*$ by $n\al_1\mapsto nQ$.  
We denote by $\C[\bar{H}^*]$ the group algebra over $\C$ of $\bar{H}^*$. 
We denote by $e^{\al}$ the element of $\C[\bar{H}^*]$ corresponding to $\al\in \bar{H}^*$. 
These $e^\al$ satisfy $e^\al e^\beta=e^{\al+\beta}$
 and $(e^\al)^{-1}=e^{-\al}$. 
In particular, $e^0=1$ is the identity element.

Now we take the power series field  $\FF=\C((P,r^*))$ 
as $\K$ and consider the semi-direct product $\C$-algebra 
$U_{q,p}(\slth)=\FF[U_q(\slth)]\otimes_{\C} \C[\bar{H}^*]$ of 
$\FF[U_q(\slth)]$ and $\C[\bH^*]$. 
We impose the following relation. For $\al\in \bH^*$, 
\bea
&&e^{\al}f(P,r^*)e^{-\al}=f(P+<\al,P>,r^*).\lb{PinvQ}
\ena
Then the  multiplication of $U_{q,p}(\slth)$ is defined by
\be
&&(f(P,r^*)a\otimes e^{\al})\cdot(g(P,r^*)b\otimes e^{\beta})=f(P,r^*)g(P+<\al,P>,r^*)ab\otimes e^{\al+\beta},\\
&&a,\ b\in \C[U_q(\slth)],\  f(P,r^*), g(P,r^*)\in \FF,\ \al, \beta\in \bH^*. 
\en
Moreover, for $f(P,r^*)\in \FF$ we regard the object $f(P+h,r^*+c)$ as the element of 
$U_{q,p}(\slth)$ in the sense of completion. We then have the following relations.
\bea
&&x^{\pm}(z)f(P+h,r^*+c)=f(P+h\mp 2,r^*+c)x^{\pm}(z),\lb{hinvXpm}\\ 
&&[f(P+h,r^*+c),a_n]=0,\qquad [f(P+h,r^*+c),d]=0.
\ena
\noi
{\it Remark.} The relation \eqref{PinvQ} is automatically satisfied, if one takes the realization 
$Q=\frac{\partial }{\partial P}$.

The following automorphism $\phi_r$ of $\FF[U_q(\slth)]$ is the key to our ``elliptic deformation''
\cite{JKOS2}.  
\be
&&c\mapsto c,\quad h\mapsto h, \quad d\mapsto d,\\
&&x^+(z)\mapsto  u^+(z,p)x^+(z),\quad x^-(z)\mapsto x^-(z)u^-(z,p),\\
&&\psi(z)\mapsto u^+(q^{c/2}z,p)\psi(z)u^-(q^{-c/2}z,p),\\
&&\varphi(z)\mapsto u^+(q^{-c/2}z,p)\varphi(z)u^-(q^{c/2}z,p).
\en
Here we set 
\bea
&&
\hspace*{-1cm}u^+(z,p)
=\exp\left(\sum_{n>0}\frac{1}{[r^*n]_q}a_{-n}(q^rz)^n\right),\lb{uplus} 
\quad u^-(z,p)=\exp\left(-\sum_{n>0}\frac{1}{[rn]_q}a_{n}(q^{-r}z)^{-n}\right).\lb{uminus}
\ena

We define the elliptic currents $E(u), F(u)$ and $K(u)$ in $U_{q,p}(\slth)[[u]]$ as follows. 
\begin{dfn}[Elliptic currents]\lb{ecurrents}
\be
\!\!\! E(u)&=&\phi_r(x^+(z))e^{2Q}z^{-\frac{P-1}{r^*}},\\
\!\!\! F(u)&=&\phi_r(x^-(z))z ^{\frac{P+ h-1}{r}},\\
\!\!\!K(z)&=&\exp\left(\sum_{n>0}\frac{[n]_q}{[2n]_q[r^*n]_q}a_{-n}(q^cz)^n\right)
\exp\left(-\sum_{n>0}\frac{[n]_q}{[2n]_q[rn]_q}a_{n}z^{-n}\right)\\
&&\qquad\qquad\times e^Qz^{-\frac{c}{4rr^*}(2P-1)+\frac{1}{2r}h}
\\
\!\!\!\hat{d}&=&d-\frac{1}{4r^*}(P-1)(P+1)+\frac{1}{4r}(P+h-1)(P+h+1), 
\en
where we set $z=q^{2u}$, $p=q^{2r}$.
\end{dfn}
From \eqref{Heisenberg} and Definition \ref{defUq}, 
we can derive the following relations. 
\begin{prop}\lb{Defrels} 
\be
&&{c}:\hbox{ {\rm central}}, 
\lb{u1}
\\
&&[h,a_{n}]=0,\quad [h,E(u)]=2E(u),\quad 
[h,F(u)]=-2 F(u), 
\lb{u2}\\
&&[\hat{d},h]=0,
\quad [\hat{d},a_{n}]= n a_{n},\quad\\ 
&&[\hat{d},E(u)]=\left(-z\frac{\partial}{\partial z}-\frac{1}{r^*}\right)
E(u) , \quad
[\hat{d},F(u)]=\left(-z\frac{\partial}{\partial z}-\frac{1}{r}\right)
F(u) ,
\lb{u3}\\
&&
[a_{n},a_{m}]=\frac{[2n]_q[c n]_q}{n}
q^{-c|n|}\delta_{n+m,0},
\label{u4}\\
&&
[a_{n},E(u)]=\frac{[2n]_q}{n}q^{-c|n|}z^n E(u),
\quad
[a_{n},F(u)]=-\frac{[2n]_q}{n} z^n F(u),
\lb{u6}\\
&& E(u)E(v)
=\frac{\left[u-v+1\right]^*}{\left[u-v-1\right]^* } 
E(v)E(u),
\quad
 F(u)F(v)
=\frac{\left[u-v-1\right]}{\left[u-v+1 \right]}
F(v)F(u),
\lb{u8}\\
&&[E(u),F(v)] =\frac{1}{q-q^{-1}}
\left(\delta\left(q^{-c}\frac{z}{w}\right)H^+(q^{c/2}w)
-\delta\left(q^{c}\frac{z}{w}\right)H^-(q^{-c/2}w)
\right),\lb{u9}
\en
where $z=q^{2u}$, $w=q^{2v}$, and we set  
\bea
&&H^\pm(z)=\kappa K\left(u\pm\frac{1}{2}(r-\frac{c}{2})+\frac{1}{2}\right)
 K\left(u\pm\frac{1}{2}(r-\frac{c}{2})-\frac{1}{2}\right),\lb{HKK}\\
&&\kappa=\lim_{z\to q^{-2}}\frac{\xi(z;p^*,q)}{\xi(z;p,q)},\qquad 
\xi(z;p,q)=\frac{(q^2z;p,q^4)_\infty(pq^2z;p,q^4)_\infty}
{(q^4z;p,q^4)_\infty(pz;p,q^4)_\infty}.\nn
\ena
\end{prop}
Moreover from \eqref{PinvQ} and \eqref{hinvXpm}, we obtain the following relations. 
\begin{prop}\lb{commKEFf}
For $f(P)\in \C((P))$, 
\be
&&K(u)f(P)=f(P+1)K(u),\quad E(u)f(P)=f(P+2)E(u),\quad [F(u),f(P)]=0,
\\
&&K(u)f(P+h)=f(P+h+1)K(u),\quad [E(u),f(P+h)]=0,\quad F(u)f(P+h)=f(P+h+2)F(u).
\en
\end{prop}
\begin{dfn}
We call a set ( $\FF[U_q(\slth)]\otimes_{\C} \C[\bH^*]$,  $\phi_r$)  
the elliptic algebra $U_{q,p}(\slth)$.  
We also denote  by $U'_{q,p}(\slth)$ the subalgebra $\FF[U'_q(\slth)]\otimes_{\C} \C[\bH^*]$ 
of $U_{q,p}(\slth)$.    
\end{dfn}

The following relations are also useful. 
\begin{prop}\lb{relecurrents2}
\be
\!\!\!\!\!\!\!\!\!\!&&K(u)K(v)=\rho(u-v)K(v)K(u),
\lb{Ucom1}\\[2mm]
\!\!\!\!\!\!\!\!\!\!&&K(u)E(v)=\frac{[u-v+\frac{1-r^*}{2}]^*}
{[u-v-\frac{1+r^*}{2}]^*}E(v)K(u),
\qquad K(u)F(v)=\frac{[u-v-\frac{1+r}{2}]}
{[u-v+\frac{1-r}{2}]}F(v)K(u),
\\
&&H^{+}(u)H^{-}(v)=\frac{[{u-v-1-\frac{c}{2}}]}{
[{u-v+1-\frac{c}{2}}]}
\frac{[{u-v+1+\frac{c}{2}}]^*}{
[{u-v-1+\frac{c}{2}}]^*}H^{-}(v)H^{+}(u),\\
&&H^{\pm}(u)H^{\pm}(v)=\frac{[{u-v-1}]}{
[{u-v+1}]}
\frac{[{u-v+1}]^*}{
[{u-v-1}]^*}H^{\pm}(v)H^{\pm}(u),
\en
where
\be
&&\rho(u)=\frac{\rho^{+*}(u)}{\rho^+(u)},\quad 
\rho^{+*}(u)=\rho^+(u)|_{r\to r^*},\lb{rho11}\\
&&\rho^+(u)=z^{\frac{1}{2r}}
\frac{\{pq^2z\}^2}{\{pz\}\{pq^4z\}}
\frac{\{z^{-1}\}\{q^4z^{-1}\}}{\{q^2z^{-1}\}^2},\quad \{z\}=(z;p,q^4)_\infty.
\nn
\en
\end{prop}
Note that the function $\rho(u)$  satisfies
\begin{eqnarray*}
&&\rho(0)=1,
\quad 
\rho(1)=\frac{[1]^*}{[1]},\quad \rho(u)\rho(-u)=1,
\quad
\rho(u)\rho(u+1)=\frac{[u+1]^*}{[u]^*}
\frac{[u]}{[u+1]}.
\end{eqnarray*}

In addition, the following formulae indicate  
a direct construction of $H^\pm(u)$ from the Drinfeld currents $\psi(z)$ and $\varphi(z)$.
\begin{prop}\lb{hpm}
\be
H^+(u)&=&\phi_r(\psi(z))e^{2Q}\left(q^{r-c/2}z\right)^{-\frac{c}{rr^*}(P-1)+\frac{h}{r}},\\
H^-(u)&=&\phi_r(\varphi(z))e^{2Q}\left(q^{-r+c/2}z\right)^{-\frac{c}{rr^*}(P-1)+\frac{h}{r}}.
\en
\end{prop}

\subsection{The $RLL$-relation for $U_{q,p}(\slth)$}

Following  \cite{JKOS2}, we summarize the results on the $L$ operator 
and the $RLL$-relation for $U_{q,p}(\slth)$.
 In Sect. \ref{hHopfdef}, we use the $L$ operator to define 
the $\h$-Hopf algebroid structure of $U_{q,p}(\slth)$. 

We first define the half currents $E^+(u), F^+(u)$ and $K^+(u)$ as follows.
\begin{dfn}[Half currents]\lb{halfcurrents}
\bea
&&
K^+(u)=K(u+\tfrac{r+1}{2}),
\label{kplush}\\
&&E^+(u)
=a^* \oint_{C^*} E(u') 
\frac{\left[u-u'+c/2-P+1\right]^*[1]^*}
{[u-u'+c/2]^*[P-1]^*}
\frac{dz'}{2\pi i z'},
\lb{Eplus}\\
&&F^+(u)
=a \oint_{C} F(u') 
\frac{\left[u-u'+P+h-1\right][1]}{[u-u'][P+h-1]}
\frac{dz'}{2\pi i z'}. 
\lb{Fplus}
\ena
Here the contours are chosen such that 
\be
C^* &:& |p^*q^c z|<|z'|<|q^cz|, 
\qquad C : |pz|<|z'|<|z|,
\lb{C}
\en
and the constants $a,a^*$ are chosen to satisfy 
\begin{eqnarray*}
{a^* a [1]^*\kappa\over q-q^{-1}} =1.
\end{eqnarray*}
\end{dfn}

The commutation relations for the elliptic 
currents in Propositions \ref{Defrels}$-$\ref{relecurrents2} yield  the following relations for the half currents.

\begin{prop}\lb{relHalfcurrents}
\bea
\!\!\!\!\!\!\!\!\!\!&&K^+(u_1)K^+(u_2)=\rho(u)K^+(u_2)K^+(u_1),
\lb{hf1}\\
\!\!\!\!\!\!\!\!\!\!&&
K^+(u_1)E^+(u_2)K^+(u_1)^{-1}=E^+(u_2)\frac{[1+u]^*}{[u]^*}
-E^+(u_1){\left[1\right]^*\over\left[P\right]^*}
\frac{[P+u]^*}{[u]^*},
\lb{hf2}\\
\!\!\!\!\!\!\!\!\!\!&&
K^+(u_1)^{-1}F^+(u_2)K^+(u_1)
=
\frac{[1+u]}{[u]}F^+(u_2)
-
{\left[1\right]\over
\left[P+h\right]}
\frac{\left[P+h-u\right]}{[u]}F^+(u_1),
\ena
\bea
\!\!\!\!\!\!\!\!\!\!&&
\frac{[1-u]^*}{[u]^*}E^+(u_1)E^+(u_2)
+
\frac{[1+u]^*}{[u]^*}E^+(u_2)E^+(u_1)
\lb{hf4}\\
\!\!\!\!\!\!\!\!\!\!&&
\qquad=
E^+(u_1)^2
{\left[1\right]^* \over
\left[P-2\right]^*}
{\left[P-2+u\right]^* \over
\left[u\right]^*}
+E^+(u_2)^2{\left[1\right]^* \over
\left[P-2\right]^*}
{\left[P-2-u\right]^* \over
\left[u\right]^*},
\nonumber\\
\!\!\!\!\!\!\!\!\!\!&&
\frac{[1+u]}{[u]}F^+(u_1)F^+(u_2)
+\frac{[1-u]}{[u]}F^+(u_2)F^+(u_1)
\lb{hf5} \\
\!\!\!\!\!\!\!\!\!\!&&=
F^+(u_1)^2
{\left[1\right] \over
\left[P+h-2\right]}
{\left[P+h-2-u\right] \over
\left[u\right]}
+F^+(u_2)^2
{\left[1\right] \over
\left[P+h-2\right]}
{\left[P+h-2+u\right] \over
\left[u\right]},
\nonumber\\
\!\!\!\!\!\!\!\!\!\!&&[E^+(u_1),F^+(u_2)]=
K^+(u_2-1)K^+(u_2)
\frac{\left[P-1-u\right]^*}{[u]^*}
\frac{[1]^*}{[P-1]^*} 
\nonumber\\
\!\!\!\!\!\!\!\!\!\!&&\qquad\qquad\qquad\quad\quad -
K^+(u_1)K^+(u_1-1)
\frac{\left[P+h-1-u\right]}{[u]}
\frac{[1]}{[P+h-1]},
\lb{hf6}
\ena
where we set $u=u_1-u_2$. 
\end{prop}

We next define the $L$-operator  $\hL^+(u)\in \End(V)\otimes\Uqp\bigl(\slth\bigr)$ 
with $V\cong\C^2$ as follows. 

\begin{dfn}[$L$-operator]\lb{Loperator}
\bea
&&\hL^+(u)=
\left(
\begin{array}{cc}
1 &F^+(u) \\
0 &1            \\
\end{array}
\right)
\left(
\begin{array}{cc}
K^+(u-1) & 0 \\
0             &K^+(u)^{-1}\\
\end{array}
\right)
\left(
\begin{array}{cc}
1            &0      \\
E^+(u) &1      \\
\end{array}
\right).
\lb{Gauss}
\ena
\end{dfn}

Then the relations in Proposition \ref{relHalfcurrents} can be combined into the following 
single $RLL$ relation. 
\begin{prop}\lb{prop:RLL1}
The $\hL^+(u)$ operator satisfies the following $RLL$ relation. 
\bea
&&\hspace*{-1cm}R^{+(12)}(u_1-u_2,P+h)
\hL^{+(1)}(u_1)\hL^{+(2)}(u_2)
=
\hL^{+(2)}(u_2)\hL^{+(1)}(u_1)
R^{+*(12)}(u_1-u_2,P),\lb{RLL3}
\ena
where  $R^{+}(u,P+h)$ and $R^{+*}(u,P)=R^{+}(u,P)|_{r\to r^*}$ 
denote the elliptic dynamical $R$ matrices given by 
\bea
&&R^+(u,s)=\rho^+(u)
\left(
\begin{array}{cccc}
1 &                  &             & \\
  &b(u,s)      &c(u,s) & \\
  &\bar{c}(u,s)&\bar{b}(u,s) & \\
  &                  &             &1 \\
\end{array}
\right)
\lb{Rmat}
\ena
with $\rho^+(u)$ in Proposition \ref{relecurrents2}, and 
\be
&&b(u,s)=
\frac{[s+1] [s-1]  }{[s]^2}
\frac{[u]}{[1+u]},
\quad 
c(u,s)=
\frac{[1]}{[s]}
\frac{[s+u]}{[1+u]},
\lb{Rmat1}\\
&&\bar{c}(u,s)=\frac{[1]}{[s]}
\frac{[s-u]}{[1+u]},
\quad \qquad\qquad\quad\; 
\bar{b}(u,s)=
\frac{[u]}{[1+u]}.
\lb{Rmat4}
\en
\end{prop}

One should note that the $c=0$ case of the $RLL$ relation \eqref{RLL3} is 
identical,  up to a gauge transformation, with the one studied in the formulation of 
 Felder's elliptic quantum group in \cite{EV1,KNR,TV}.

\subsection{Connection to the Quasi-Hopf Algebra $\Bqla(\slth)$}

It is worth to remark a connection of  $U_{q,p}(\slth)$  to the quasi-Hopf algebra  $\Bqla(\slth)$.

Let us define a new $L$ operator $L^+(u,P)$ by
\bea
L^+(u,P)
&=&\hL^+(u)
\left(\begin{array}{cc}
e^{-Q} &0\\
0 &e^{Q}\\
\end{array}\right).\lb{dynamicalL}
\ena
Then from Definitions \ref{ecurrents} and \ref{Loperator}, one finds
 that  $L^+(u,P)$ is independent of ${Q}$. 
We hence regard $L^+(u,P)$ as the operator in $\FF[\uq]$ having $P$ as a parameter. 
Substituting \eqref{dynamicalL} into  \eqref{RLL3}, we  obtain the following statement.
\begin{prop}\lb{prop:RLL2} 
The operator $L^+(u,P)$ satisfies the following dynamical $RLL$ relation. 
\bea
&&R^{+(12)}(u_1-u_2,P+h)
L^{+(1)}(u_1,P)L^{+(2)}(u_2,P+h^{(1)})
\lb{RLL4}\\
&&\qquad =
L^{+(2)}(u_2,P)L^{+(1)}(u_1,P+h^{(2)})
R^{+*(12)}(u_1-u_2,P).
\nn
\ena
\end{prop}
This $RLL$ relation is identified with the one for 
the quasi-Hopf algebra $\Bqla(\slth)$ under the parametrization 
$\la=(r^*+2)\Lambda_0+(P+1)\bar{\Lambda}_1$ \cite{JKOS, JKOS2}, where $\Lambda_0$ and 
$\Lambda_0+\bar{\Lambda}_1$ denote the fundamental weights of $\widehat{\gsl}(2,\C)$. 
This is due to the fact that under this $\la$   
 the vector representation of the universal dynamical $R$ matrix $\cR^+(\la)$ of 
$\Bqla(\slth)$  yields the elliptic dynamical $R$ matrix $R^{+*}(u,P)$
\cite{JKOS2,Konno06}.   
Furthermore we have the isomorphism $\Bqla(\slth)\cong \FF[\uq]$ as an associative algebra. 
Combining these facts, we obtain  the isomorphism 
$U_{q,p}(\slth)\cong \Bqla(\slth)\otimes_{\C} \C[\bar{H}^*]$ 
with $\la=(r^*+2)\Lambda_0+(P+1)\bar{\Lambda}_1$ as a semi-direct product 
algebra. 

Note also that the $c=0$ case of \eqref{RLL4} is identical to the one used in \cite{Felder,FV} to 
define Felder's elliptic quantum group in its original form. 
\section{$H$-Hopf Algebroid}
In this section, we introduce an $\h$-Hopf algebroid structure into  
the elliptic algebra $U_{q,p}(\slth)$ and formulate $U_{q,p}(\slth)$ as an elliptic 
quantum group.

\subsection{Definition of the $\h$-Hopf Algebroid}
Let us recall some basic facts on the $\h$-Hopf algebroid following the works of  
Etingof and Varchenko\cite{EV1,EV2} and of Koelink and Rosengren \cite{KR}. 

Let $A$ be a complex associative algebra, $\h$ be a finite dimensional commutative subalgebra of 
$A$, and $M_{\h^*}$ be the 
field of meromorphic functions on $\h^*$ the dual space of $\h$. 

\begin{dfn}[$\h$-algebra] 
An $\h$-algebra is a complex associative algebra $A$ with 1, which is bigraded over 
$\h^*$, $\ds{A=\bigoplus_{\alpha,\beta\in \h^*} A_{\al\beta}}$, and equipped with two 
algebra embeddings $\mu_l, \mu_r : M_{\h^*}\to A_{00}$ (the left and right moment maps), such that 
\be
\mu_l(\hf)a=a \mu_l(T_\al \hf), \quad \mu_r(\hf)a=a \mu_r(T_\beta \hf), \qquad 
a\in A_{\al\beta},\ \hf\in M_{\h^*},
\en
where $T_\al$ denotes the automorphism $(T_\al \hf)(\la)=\hf(\la+\al)$ of $M_{\h^*}$.
\end{dfn}
\begin{dfn}[$\h$-algebra homomorphism] 
An $\h$-algebra homomorphism is an algebra homomorphism $\pi:A\to B$ between two $\h$-algebras $A$ and $B$ preserving the bigrading and the moment maps, i.e. $\pi(A_{\al\beta})\subseteq B_{\al\hb}$ for all $\ha,\hb\in \h^*$ and $\pi(\mu^A_l(\hf))=\mu^B_l(\hf), \pi(\mu^A_r(\hf))=\mu^B_r(\hf)$. 
\end{dfn}

Let $A$ and $B$ be two $\h$-algebras. The tensor product $A {\widetilde{\otimes}}B$ is the $\h^*$-bigraded vector space with 
\be
 (A {\widetilde{\otimes}}B)_{\al\beta}=\bigoplus_{\gamma\in\h^*} (A_{\al\gamma}\otimes_{M_{\h^*}}B_{\gamma\beta}),
\en
where $\otimes_{M_{\h^*}}$ denotes the usual tensor product 
modulo the following 
relation.
\bea
\mu_r^A(\hf) a\otimes b=a\otimes\mu_l^B(\hf) b, \qquad a\in A, 
b\in B, \hf\in M_{\h^*}.\lb{AtotB}
\ena
The tensor product $A {\widetilde{\otimes}}B$ is again an $\h$-algebra with the multiplication $(a\otimes b)(c\otimes d)=ac\otimes bd$ and the moment maps 
\be
\mu_l^{A {\widetilde{\otimes}}B} =\mu_l^A\otimes 1,\qquad \mu_r^{A {\widetilde{\otimes}}B} =1\otimes \mu_r^B.
\en

Let $\cD$ be the algebra of automorphisms $M_{\h^*}\to M_{\h^*}$ 
\be
\cD&=&\{\ \sum_i \hf_i T_{\beta_i}\ |\ \hf_i\in M_{\h^*},\ \beta_i\in \h^*\ \}.
\en
Equipped  with the bigrading 
 $\cD_{\al\al}=\{\  \hf T_{-\al}\ |\ \hf\in M_{\h^*},\ \al\in \h^*\ \}$, 
 $\cD_{\al\beta}=0\ (\al\not=\beta)$ 
 and the moment maps $\mu^{\cD}_l, \mu^{\cD}_r : M_{\h^*}\to \cD_{00}$ 
 defined by 
$\mu^{\cD}_l(\hf)=\mu^{\cD}_r(\hf)=\hf T_0$, $\cD$ is an $\h$-algebra.
 For any $\h$-algebra $A$, we have the canonical 
isomorphism as an $\h$-algebra 
\bea&&
A\cong A\tot \cD\cong  \cD\tot A \lb{Diso}
\ena
by $a\cong a\tot T_{-\beta}\cong T_{-\al}\tot a$ for all $a\in A_{\al\beta}$.

\begin{dfn}[$\h$-bialgebroid]
An $\h$-bialgebroid is an $\h$-algebra $A$ equipped with two $\h$-algebra homomorphisms 
$\Delta:A\to A{{\tot}}A$ (the comultiplication) and $\vep : A\to \cD$ (the counit) such that 
\be
&&(\Delta \tot \id)\circ \Delta=(\id \tot \Delta)\circ \Delta,\\
&&(\vep \tot \id)\circ\Delta =\id =(\id \tot \vep)\circ \Delta,
\en
under the identification \eqref{Diso}.
\end{dfn}
 
\begin{dfn}[$\h$-Hopf algebroid]\lb{defS}
An $\h$-Hopf algebroid is an $\h$-bialgebroid $A$ equipped with a $\C$-linear map $S : A\to A$ (the antipode), such that 
\be
&&S(\mu_r(\hf)a)=S(a)\mu_l(\hf),\quad S(a\mu_l(\hf))=\mu_r(\hf)S(a),\quad \forall a\in A, \hf\in M_{\h^*},\\
&&m\circ (\id \tot S)\circ\Delta(a)=\mu_l(\vep(a)1),\quad \forall a\in A,\\
&&m\circ (S\tot\id  )\circ\Delta(a)=\mu_r(T_{\al}(\vep(a)1)),\quad \forall a\in A_{\al\beta},
\en
where $m : A{{\tot}} A \to A$ denotes the multiplication and $\vep(a)1$ is the result of applying the difference operator $\vep(a)$ to the constant function $1\in M_{\h^*}$.
\end{dfn}

\noi
{\it Remark.}\cite{KR} Definition \ref{defS} yields that the antipode of an $\h$-Hopf algebroid 
 uniquely exists and gives the algebra antihomomorphism.

The $\h$-algebra $\cD$ is an $\h$-Hopf algebroid with 
$\Delta_\cD : \cD\to \cD\tot \cD,\ \vep_\cD: \cD \to \cD,\ 
S_\cD : \cD \to \cD$ defined by 
\be
&&\Delta_\cD(\hf T_{-\al})=\hf T_{-\al} \tot T_{-\al},\\
&&\vep_\cD=\id,
\qquad  S_\cD(\hf T_{-\al})=T_{\al}\hf=(T_{\al}\hf)T_{\al}.
\en

\subsection{$H$-Hopf Algebroid Structure on $U_{q,p}(\slth)$}\lb{hHopfdef}
Now let us consider the elliptic algebra $U_{q,p}(\slth)$.  
Using the isomorphism  $\phi:\cQ\to \bar{\h}^*$, we define 
the $\bar{\h}^*$-bigrading structure of $U_{q,p}=U_{q,p}(\slth)$ as follows.  
\bea
&&{U}_{q,p}=\bigoplus_{\al, \beta \in \bar{\h}^*} (U_{q,p})_{\al \beta},\nn\\
&&(U_{q,p})_{\al\beta }=\left\{\ x\in U_{q,p}\ \left|\ 
\mmatrix{q^{h}xq^{-h}=q^{<\phi^{-1}({\al-\beta}),h>}x \cr
q^{P}xq^{-P}=q^{<\beta,P>}x \cr}\ \right.\right\}.\lb{bigrading}
\ena
Noting $<\phi^{-1}({\al}),h>=<{\al},P>$, we have 
\bea
q^{P+h}xq^{-(P+h)}=q^{<\al,P>}x \lb{gradePh}
\ena 
for $x\in (U_{q,p})_{\al\beta }$. 

\noi
{\it Remark.} The quantum affine algebra $U_q=\FF[U_q(\slth)]$ has 
the following natural grading over $\bH^*$. 
\be
U_q&=&\bigoplus_{{\al}\in \bH^*} (U_q)_{{\al}},\qquad 
(U_q)_{{\al}}=\{x\in U_q\ |\ q^hxq^{-h}=q^{<\phi^{-1}({\al}),h>}x\ 
\}.
\en
We then have 
\be
&&(U_{q,p})_{\ha\hb}=(U_q)_{{\al-\beta}}\otimes_{\C} \C e^{-\beta}.
\en

 Next let us regard the elements 
$\widehat{f}=f(P,r^*)\in \FF$ as  meromorphic functions on $\h^*$ by 
\be
&&\widehat{f}(\mu)=f(<\mu,P>,<\mu,r^*>)\quad \mu\in {H}^*
\en
and consider  the field of  meromorphic functions  $M_{{\h}^*}$ on ${\h}^*$ 
\be
&&{M}_{{\h}^*}=\left\{ \widehat{f}:{\h}^*\to \C\ \left|\ 
\widehat{f}=f(P,r^*)\in \FF\right.\right\}.
\en 
We define two embeddings (the left and right moment maps) 
$\mu_l, \mu_r : {M}_{{\h}^*}\to (U_{q,p})_{00}$  by
\bea
\mu_l(\widehat{f})=f(P+h,{r^*+c}),\qquad 
\mu_r(\widehat{f})=f(P,{r^*}).\lb{mmUqp}
\ena
From \eqref{PinvQ} and \eqref{hinvXpm}, one can verify the following. 
\begin{prop}\lb{commMu}
For $x\in (U_{q,p})_{\al\beta}$, we have 
\be
&&\mu_l(\widehat{f})x=f(P+h,r^*+c)x=x f(P+h+<\al,P>,r^*+c)=x \mu_l({T}_{\al}
\widehat{f}),\lb{mlUqp}\\
&&\mu_r(\widehat{f})x=f(P,r^*)x=x f(P+<\beta,P>,r^*)=x \mu_r({T}_{\beta}\widehat{f}),\qquad \lb{mrUqp}
\en
where we regard ${T}_{\al}=e^\al\in \C[\bH^*]$ as 
the shift operator $M_{{\h}^*}\to M_{{\h}^*}$ 
\be
({T}_{\al}\widehat{f})=e^{\al}f(P,r^*)e^{-\al}={f}(P+<\al,P>,r^*).
\en
\end{prop}
Hereafter we abbreviate 
$f(P+h,{r^*+c})$ and $f(P,{r^*})$ as $f(P+h)$ and
 $f^*(P)$, respectively.  

An important example of the elements in 
${M}_{\h^*}$ is the 
elliptic dynamical $R$ matrix elements 
$(\widehat{R}^+_u)_{\vep_1\vep_2}^{\vep_1'\vep_2'}\equiv R^{+*}(u,P)_{\vep_1\vep_2}^{\vep_1'\vep_2'}$ in \eqref{Rmat}, where $\vep_i, \vep'_i=+, -\ (i=1,2)$. 
We then have 
\bea
&&\mu_l((\widehat{R}^+_u)_{\vep_1\vep_2}^{\vep_1'\vep_2'})=R^+(u,P+h)_{\vep_1\vep_2}^{\vep_1'\vep_2'},\quad  
\mu_r((\widehat{R}^+_u)_{\vep_1\vep_2}^{\vep_1'\vep_2'})
=R^{+*}(u,P)_{\vep_1\vep_2}^{\vep_1'\vep_2'}\lb{RmatMH}
\ena
in the abbreviate notation.

 Equipped with the 
bigrading structure \eqref{bigrading} and two moment maps \eqref{mmUqp}, 
the elliptic algebra $U_{q,p}(\slth)$ is an $\h$-algebra.

We also consider the $\h$-algebra of the shift operators 
\be
&&\cD=\{\ \sum_i \widehat{f}_i{T}_{\ha_i}\ |\ \widehat{f}_i \in {M}_{\h^*}, 
\ha_i\in \bH^*\ \},\\
&&\cD_{\ha\ha}=\{\ \widehat{f}{T}_{-\ha}\ \},\quad \cD_{\ha{\beta}}=0\ 
(\al\not=\beta),\\
&&\mu_l^{\cD}(\widehat{f})=
\mu_r^{\cD}(\widehat{f})=\widehat{f}{T}_0 \qquad \widehat{f}\in {M}_{\h^*}.
\en
Then we have the $\h$-algebra isomorphism $ U_{q,p}\cong U_{q,p}\tot\cD\cong \cD\tot U_{q,p}$. 

Now let us consider the $\h$-Hopf algebroid structure on $U_{q,p}$. 
It is conveniently given by the $L$ operator $\hL^+(u)$. 
We shall write the entries of $\widehat{L}^+(u)$ as 
\bea
&&\hLp(u)=\left(\mmatrix{\hLp_{++}(u)&\hLp_{+-}(u)\cr
\hLp_{-+}(u)&\hLp_{--}(u)\cr
}\right).\lb{Lmatelements}
\ena
According to the Gau{\ss} decomposition \eqref{Gauss}, we have
\bea
&&\widehat{L}^+_{++}(u)=K^+(u-1)+F^+(u)K^+(u)^{-1}E^+(u),\quad
\widehat{L}^+_{+-}(u)=F^+(u)K^+(u)^{-1},\lb{Lelements}
\\
&&\widehat{L}^+_{-+}(u)=K^+(u)^{-1}E^+(u),\qquad\qquad\qquad\qquad\quad\;\;
\widehat{L}^+_{--}(u)=K^+(u)^{-1}.\nn
\ena
 One finds 
\bea
&&\hL^+_{\vep_1\vep_2}(u)\in (U_{q,p})_{-{\vep_1Q},-\vep_2Q}.
\lb{grL1}
\ena
It is also easy to check 
\bea
f(P+h)\hL^+_{\vep_1\vep_2}(u)&=&
\hL^+_{\vep_1\vep_2}(u)f(P+h-\vep_1),\nn\\
f^*(P)\hL^+_{\vep_1\vep_2}(u)&=&\hL^+_{\vep_1\vep_2}(u)
f^*(P-\vep_2).\lb{shifts}
\ena

We define two $\h$-algebra homomorphisms, the co-unit $\vep : U_{q,p}\to \cD$ and the co-multiplication $\Delta : U_{q,p}\to U_{q,p} \widetilde{\otimes}U_{q,p}$ by
\bea
&&\vep(\hL^+_{\vep_1\vep_2}(u))=\delta_{\vep_1,\vep_2}{T}_{\vep_2 },
\quad \vep(e^Q)=e^Q,\lb{counitUqp}\\
&&\vep(\mu_l({\hf}))= \vep(\mu_r(\hf))=\widehat{f}T_0, \lb{counitf}\\
&&\Delta(\hL^+_{\vep_1\vep_2}(u))=\sum_{\vep'}\hL^+_{\vep_1\vep'}(u)\widetilde{\otimes}
\hL^+_{\vep'\vep_2}(u),\lb{coproUqp}\\
&&\Delta(e^{Q})=e^{Q}\tot e^{Q},\\
&&\Delta(\mu_l(\hf))=\mu_l(\hf)\widetilde{\otimes} 1,\quad \Delta(\mu_r(\hf))=1\widetilde{\otimes} \mu_r(\hf).\lb{coprof}
\ena
In fact, one can check that $\Delta$ preserves the relation \eqref{RLL3}. 
Noting \eqref{RmatMH} 
and the formula obtained from \eqref{AtotB}
\bea
&&f^*(u,P)a\tot b=a\tot f(u,P+h)b\qquad a,b\in U_{q,p},\lb{fstotf}
\ena
 we have
\be
\Delta(LHS)&=&\sum_{\vep_1',\vep_2'}
\Delta(R^+(u,P+h)_{\vep_1''\vep_2''}^{\vep_1'\vep_2'})
\Delta(\hL^+_{\vep_1'\vep_1}(u_1))\Delta(\hL^+_{\vep_2'\vep_2}(u_2))\\
&=&\sum_{\vep_1',\vep_2'\atop \vep, \vep'}
R^+(u,P+h)_{\vep_1''\vep_2''}^{\vep_1'\vep_2'}\hL^+_{\vep_1'\vep}(u_1)\hL^+_{\vep_2'\vep'}(u_2)\widetilde{\otimes} \hL^+_{\vep\vep_1}(u_1)\hL^+_{\vep'\vep_2}(u_2)\\
&=&\sum_{\vep_1',\vep_2'\atop \vep, \vep'}
\hL^+_{\vep_2''\vep'_2}(u_2)\hL^+_{\vep_1''\vep_1'}(u_1)
R^{+*}(u,P)_{\vep_1'\vep_2'}^{\vep\vep'}\widetilde{\otimes} \hL^+_{\vep\vep_1}(u_1)\hL^+_{\vep'\vep_2}(u_2)\\
&=&\sum_{\vep_1',\vep_2'\atop \vep, \vep'}R^{+*}(u,P)_{\vep_1'\vep_2'}^{\vep\vep'}
\hL^+_{\vep_2''\vep'_2}(u_2)\hL^+_{\vep_1''\vep_1'}(u_1)
\widetilde{\otimes} \hL^+_{\vep\vep_1}(u_1)\hL^+_{\vep'\vep_2}(u_2)\\
&=&\sum_{\vep_1',\vep_2'\atop \vep, \vep'}
\hL^+_{\vep_2''\vep'_2}(u_2)\hL^+_{\vep_1''\vep_1'}(u_1)
\widetilde{\otimes} R^{+}(u,P+h)_{\vep_1'\vep_2'}^{\vep\vep'}\hL^+_{\vep\vep_1}(u_1)\hL^+_{\vep'\vep_2}(u_2)\\
&=&\sum_{\vep_1',\vep_2'\atop \vep, \vep'}
\hL^+_{\vep_2''\vep'_2}(u_2)\hL^+_{\vep_1''\vep_1'}(u_1)
\widetilde{\otimes} \hL^+_{\vep_2'\vep'}(u_2)\hL^+_{\vep_1'\vep}(u_1)
R^{+*}(u,P)_{\vep\vep'}^{\vep_1\vep_2}\\
&=&\Delta(RHS).
\en
In the fourth line, we used the property
\be
&&R^{+*}(u,P+\vep_1'+\vep_2')_{\vep_1'\vep_2'}^{\vep_1\vep_2}=R^{+*}(u,P)_{\vep_1'\vep_2'}^{\vep_1\vep_2}.
\en

\begin{lem}\lb{counitcopro}
The maps $\vep$ and $\Delta$ satisfy
\bea
&&(\Delta\tot \id)\circ \Delta=(\id \tot \Delta)\circ \Delta,\lb{coaso}\\
&&(\vep \tot \id)\circ\Delta =\id =(\id \tot \vep)\circ \Delta.\lb{vepDelta}
\ena
\end{lem}
\noi
{\it Proof.} Straight forward. 
\qed

We also have the following formulae. 
\begin{prop}
\bea
&&\vep(q^{h})=\vep(q^{{c}})=T_0,\lb{ehc}\\
&&\Delta(q^{ h})=q^{ h}\tot q^{ h},\quad 
\Delta(q^{{c}})=q^{{c}}\tot q^{{c}},\lb{Dhc}\\
&&\Delta\left(\frac{f(P,r^*)}{f(P+h,r^*+c)}\right)=\frac{f(P,r^*)}{f(P+h,r^*+c)}\tot
\frac{f(P,r^*)}{f(P+h,r^*+c)}.\lb{fsf}
\ena
\end{prop}
\noi
{\it Proof.} \eqref{ehc} follows from \eqref{mmUqp} and  \eqref{counitf}, whereas \eqref{Dhc} follows from \eqref{mmUqp}, \eqref{coprof} and \eqref{AtotB}. For example, 
\be
&&\Delta(q^h)=\Delta(q^{P+h}q^{-P})=\Delta(q^{P+h})\Delta(q^{-P})=q^{P+h}\tot q^{-P}=q^{h}\tot q^{h}.
\en
To show \eqref{fsf} we use \eqref{coprof} and \eqref{AtotB} as
\be
{\rm LHS}=\Delta(\mu_r(\widehat{f}))\Delta(\mu_l(\widehat{f})^{-1})=
\mu_l(\widehat{f})^{-1}\tot\mu_r(\widehat{f})=\frac{f(P,r^*)}{f(P+h,r^*+c)}\frac{1}{f(P,r^*)}\tot
f(P,r^*)={\rm RHS}.
\en
\qed

We next define an algebra antihomomorphism (the antipode) $S : U_{q,p}\to U_{q,p}$ by
\be
&&S(\hL^+_{++}(u))=\hL^+_{--}(u-1), \quad S(\hL^+_{+-}(u))=-\frac{[P+h+1]}{[P+h]}
\hL^+_{+-}(u-1), \quad \\
&&S(\hL^+_{-+}(u))=-\frac{[P]^*}{[P+1]^*}\hL^+_{-+}(u-1),\quad S(\hL^+_{--}(u))=\frac{[P+h+1][P]^*}{[P+h][P+1]^*}\hL^+_{++}(u-1),\\
&&S(e^Q)=e^{-Q},\quad S(\mu_r(\hf))=\mu_l(\hf),\quad S(\mu_l(\hf))=\mu_r(\hf).
\en
Note that $S$ preserves the $RLL$ relation \eqref{RLL3}. To show this, we use 
 the relations in Proposition \ref{relHalfcurrents}. Furthermore we have the 
following Lemma.
\begin{lem}\lb{antipode}
The map $S$ satisfies 
\be
&&m\circ (\id \otimes S)\circ\Delta(x)=\mu_l(\vep(x)1),\quad \forall x\in U_{q,p},\\
&&m\circ (S\otimes\id  )\circ\Delta(x)=\mu_r(T_{\ha}(\vep(x)1)),\quad \forall x\in (U_{q,p})_{\ha \hb}.
\en
\end{lem}

\noi
{\it Proof.} We prove the first relation for $x=\hL_{++}(u)$. The other is similar.
Using \eqref{Lelements} and \eqref{shifts},
\be
LHS&=&\hL_{++}(u)\hL_{--}(u-1)-\hL_{+-}(u)\hL_{-+}(u-1)\frac{[P-1]^*}{[P]^*}\\
&=&(K^+(u-1)+F^+(u)K^+(u)^{-1}E^+(u))K^+(u-1)^{-1}\\
&&\qquad -F^+(u)K^+(u)^{-1}K^+(u-1)^{-1}E^+(u-1)\frac{[P-1]^*}{[P]^*}\\
&=&1=\mu_l(\vep(\hL_{++}(u)1)).
\en
In the last line, we used the relation \eqref{hf2} with the 
replacement $u_1\mapsto u-1, u_2\mapsto u$ and $u\mapsto -1$. 
\qed

From Lemmas \ref{counitcopro} and \ref{antipode}, we have 
\begin{thm}
The $\h$-algebra $U_{q,p}(\slth)$ equipped with $(\Delta,\vep,S)$ is an $\h$-Hopf algebroid. 
\end{thm}

\begin{dfn}
We call the $\h$-Hopf algebroid $(U_{q,p}(\slth),\h,{M}_{\h^*},\mu_l,\mu_r,\Delta,\vep,S)$ the \\ elliptic quantum group $U_{q,p}(\slth)$. 
\end{dfn}

We also use 
 the following comultiplication formulae for the half currents.
\begin{prop}\lb{coprohalfc}
\be
&&\Delta(K^+(u))=K^+(u)\tot K^+(u)+\sum_{j=1}^\infty (-)^jE^+(u)^jK^+(u)\tot 
K^+(u)F^+(u)^j,\\
&&\Delta(E^+(u))=1\tot E^+(u)+E^+(u)\tot K^+(u)K^+(u-1)\\
&&\qquad\qquad\qquad\qquad+\sum_{j=1}^\infty(-)^jE^+(u)^{j+1}\tot K^+(u)F^+(u)^jK^+(u-1),\\
&&\Delta(F^+(u))=F^+(u)\tot 1+ K^+(u-1)K^+(u)\tot F^+(u)\\
&&\qquad\qquad\qquad\qquad+\sum_{j=1}^\infty(-)^jK^+(u-1)E^+(u)^jK^+(u)\tot F^+(u)^{j+1},\\
\Delta(H^\pm(u))&=&H^\pm(u)\tot H^\pm(u)\\
&+& \sum_{j=1}^\infty(-)^j\left\{\kappa K^+(u+C_\pm)E^+(u+C_\pm-1)^j
K^+(u+C_\pm-1)\tot H^+(u) F^+(u+C_\pm-1)^j \right.\\
&&\left.\qquad\qquad\qquad +E^+(u+C_\pm)^j
H^+(u)\tot \kappa K^+(u+C_\pm) F^+(u+C_\pm)^jK^+(u+C_\pm-1)
\right\}\\
&+& \sum_{i,j=1}^\infty(-)^{i+j} \kappa E^+(u+C_\pm)^iK^+(u+C_\pm)
E^+(u+C_\pm-1)^jK^+(u+C_\pm-1)\\
&&\qquad\qquad\qquad \tot \kappa 
K^+(u+C_\pm) F^+(u+C_\pm)^iK^+(u+C_\pm-1) F^+(u+C_\pm-1)^j,
\en
where $C_\pm=-\frac{r}{2}\pm(\frac{r}{2}-\frac{c}{4})$.
\end{prop}
\noi
{\it Proof.} Use \eqref{coproUqp}, \eqref{Lelements} and \eqref{HKK} as well as 
$\Delta(\kappa)=\kappa\tot \kappa$ obtained from \eqref{fsf}. \qed

\section{Finite-Dimensional Representations}

In this section, we discuss representations of 
the elliptic algebra $U'_{q,p}=U'_{q.p}(\slth)$. 
The main results are the criterion for the finiteness of irreducible representations 
Theorem \ref{eDriPoly} and the submodule structure of the tensor product of 
two evaluation representations Theorem  \ref{singvec}$-$\ref{submodule}. 

For brevity, we denote the entries of 
$\hL^+(u)$ by 
\be
&&\hL^+(u)=\mat{\al(u)&\beta(u)\cr
                 \gamma(u)&\delta(u)\cr}.
\en

\subsection{Dynamical Representations}
We introduce the concept of dynamical 
representation, i.e. representation as $\h$-algebras\cite{EV1,EV2,KR}. 
We follows the definition given in \cite{KR}. 
We then give a construction of dynamical representations of $U'_{q,p}$. 

 Let us consider a vector space $\hV$ over $\FF$, which is  
$\bar{\hh}$-diagonalizable,
\be
&&\hV=\bigoplus_{\mu\in \bar{\hh}^*}\hV_{\mu},\qquad 
\hV_{\mu}=\{ v\in V\ |\ q^{\bar{h}}v=q^\mu v\quad (\bar{h}\in \bar{\hh})\}.
\en
Let us define the $\h$-algebra $\cD_{\h,\hV}$ of the $\C$-linear operators on $\hV$ by
\be
&&\cD_{\h,\hV}=\bigoplus_{\al,\beta\in \bH^*}(\cD_{\h,\hV})_{\al\beta},\\
&&\hspace*{-10mm}(\cD_{\h,\hV})_{\al\beta}=
\left\{\ X\in \End_{\C}\hV\ \left|\ \mmatrix{
X(f^*(P)v)=f^*(P-<\beta,P>)X(v),\cr 
X(\hV_\mu)\subseteq \hV_{\mu+\phi^{-1}(\al-\beta)},\ v\in 
\hV,\ f^*(P)\in \FF\cr}  \right.\right\},\\
&&\mu_l^{\cD_{H,\hV}}(\widehat{f})v=f(P+\mu)v,\quad 
\mu_r^{\cD_{H,\hV}}(\widehat{f})v=f^*(P)v,\qquad \widehat{f}\in {M}_{\h^*}
\en
for $v\in \hV_{\mu}$. We follow the abbreviation mentioned below Proposition \ref{commMu}. 
\begin{dfn}[Dynamical representation]
A dynamical representation of $U_{q,p}'$ on $\hV$ is 
 an $\h$-algebra homomorphism $\widehat{\pi}: U_{q,p}' 
 \to \cD_{\h,\hV}$. The dimension of the dynamical representation 
 $(\widehat{\pi},\hV)$ is  $\dim_{\FF}\hV$.  
\end{dfn}

Let 
 $(\widehat{\pi}_V, \hV), (\widehat{\pi}_W, \hW)$ be two 
dynamical representations of $U_{q,p}'$. We define the tensor product 
$\hV\tot \hW$ by
\be
&&\hV\tot \hW=\bigoplus_{\al\in \bar{\hh}^*}(\hV\tot \hW)_{\al},\quad (\hV \tot \hW)_{\al}=\bigoplus_{\beta\in  \bar{\hh}^*}\hV_{\beta}\otimes_{M_{\h^*}}\hW_{\al-\beta}, 
\en
where
$\otimes_{M_{H^*}}$ denotes the usual tensor product modulo the relation
\bea
&&f^*(P)v\otimes w=v\otimes f(P+\nu)w\lb{VtotW}
\ena
for $w\in \hW_\nu$. 
The action of the scalar $f^*(P)\in \FF$
 on the tensor space $\hV\tot \hW$ is defined as follows. 
\be
&&f^*(P).(v\tot w)=\Delta(\mu_r(\widehat{f}))(v\tot w)=v\tot f^*(P)w.
\en
We have a natural $\h$-algebra embedding $\theta_{VW}: 
\cD_{\h,\hV}\tot \cD_{\h,\hW}\to \cD_{\h,\hV\tot \hW}$ by 
$X_{\hV}\tot X_{\hW}\in 
(\cD_{\h,\hV})_{\al\gamma}\otimes_{M_{H^*}} (\cD_{\h,\hW})_{\gamma\beta} \mapsto X_{\hV}\tot 
X_{\hW}\in (\cD_{\h,\hV\tot \hW})_{\al\beta}$.   Hence 
 $\theta_{VW}\circ(\widehat{\pi}_V\otimes \widehat{\pi}_W)\circ \Delta : 
U'_{q,p}\to \cD_{\h,\hV\tot \hW}$ gives a dynamical representation 
of $U'_{q,p}$ on $\hV\tot \hW$.

Now let us consider a construction of dynamical representations of $U'_{q,p}$. 
Let $V$ be an $\bar{\hh}$-diagonalizable vector space over $\FF$.  
Let $V_Q$ be a vector space over $\C$, 
on which an action of $e^{ Q}$ is defined appropriately.  
Two important examples of $V_Q$ are  
$ V_Q=\C 1$ and $V_Q=\oplus_{n\in \Z}\C e^{nQ} $, where 
$1$ denotes the vacuum state satisfying $e^{Q}.1=1$. 
Let us consider the vector space $\hV= V\otimes_{\C} V_Q$, on which 
the actions of $f^*(P)\in \FF$ and $e^{ Q}$ are defined as follows. 
\be
&&f^*(P).(v\otimes \xi)=f^*(P)v\otimes \xi,\\
&&e^{Q}.(f^*(P)v\otimes \xi)=f^*(P+1)v\otimes e^{Q}\xi
\en
for $f^*(P)v\otimes \xi\in V\otimes V_Q$. 
The following theorem shows a construction of dynamical representations.   
\begin{thm}\lb{dynamicalrep}
Let $V, V_Q$ and $\hV$ be as in the above. 
Let $(\pi_V:\FF[U'_q]\to \End_{\FF}V, V)$ 
 be a representation of $\FF[U'_q(\slth)]$.  
Define a map $\widehat{\pi}_V=\pi_V \otimes \id: U'_{q,p}=\FF[U'_{q}]\otimes_{\C} \C[\bH^*]
\to \End_{\C}\hV$ by
\be
\!\!\! \widehat{\pi}_V( E(u))&=&\pi_V(\phi_r(x^+(z)))e^{2Q}z^{-\frac{P-1}{r^*}},\\
\!\!\! \widehat{\pi}_V(F(u))&=&\pi_V(\phi_r(x^-(z)))z ^{\frac{P+ \pi_V(h)-1}{r}},\\
\!\!\! \widehat{\pi}_V(K(u))&=&\exp\left(\sum_{n>0}\frac{[n]_q}{[2n]_q[r^*n]_q}
\pi_V(a_{-n})(q^cz)^n\right)
\exp\left(-\sum_{n>0}\frac{[n]_q}{[2n]_q[rn]_q}\pi_V(a_{n})z^{-n}\right)\\
&&\qquad\qquad\times e^Q z^{-\frac{c}{4rr^*}(2P-1)+\frac{1}{2r}\pi_V(h)}. 
\en 
Then $(\widehat{\pi}_V,\hV)$ is a dynamical representation of 
$U'_{q,p}$ on $\hV$. 
\end{thm}
Through this paper we consider the dynamical representations obtained in this way.

\subsection{Pseudo-highest Weight Representations}
We define the concept of pseudo-highest weight representations and 
write down  some basic results on them. Most of them are parallel to 
the trigonometric \cite{CP} and the rational \cite{CPY} cases.

We begin by stating an analogue of the  Poincar\'{e}-Birkhoff-Witt theorem for $U'_q$. 
\begin{dfn}
Let ${\cal H}$ (resp. ${\cal N}_\pm$) be the subalgebras of $\FF[U'_{q,p}(\slth)]$ 
generated by 
${c}, {h}$ and $a_k\ (k\in \Z_{\not=0})$ (resp. by $x^\pm_n\ (n\in \Z)$). 
\end{dfn}
From Proposition 3.1 in \cite{CP} and 
a standard normal ordering procedure on the Heisenberg algebra, we have the following. 
\begin{thm}\lb{PBW}
\be
&&U'_{q,p}=({\cal N}_-\otimes{\cal H}\otimes{\cal N}_+)\otimes \C[\bH^*].
\en
Here the last $\otimes$ should be understood as the semi-direct product. 
\end{thm}

The following indicates a characteristic feature of the
 finite-dimensional irreducible dynamical representation of $U'_{q.p}$.
\begin{thm}\lb{FDI} 
Every finite-dimensional irreducible  dynamical 
 representation $(\widehat{\pi}_V,\hV=V\otimes V_Q)$ of $U'_{q,p}$ 
  contains a non-zero vector  of the form 
$\widehat{\Omega}=\Omega\otimes 1,\ \Omega\in V$ such that
\be
1)&& x^+_n.\hOmega=0\qquad  \forall n\in \Z,  \\
2)&& {\hOmega}\ \mbox{is a simultaneous eigenvector for 
the elements of ${\cal H}$}, \\
3)&& e^{ Q}.\hOmega=\hOmega,\\
4)&& \hV=U_{q,p}'.\widehat{\Omega}.
\en
Furthermore $q^c$ acts as $1$ or $-1$ on $\hV$. 
\end{thm}
\noi
{\it Proof.} 
Note that for each $k\in \Z$, $\C\{ x_k^+, x^-_{-k}, q^hq^{kc}\}\cong U_q(\slt)$ 
is a subalgebra of $U'_{q,p}(\slth)$.  Then, concerning the action of the  
$\FF[U'_q(\slth)]$ part, the existence of 
a vector $\hOmega'=\Omega\otimes \xi\in \hV=V\otimes V_Q$ satisfying $1)$ and $2)$ 
follows from Proposition 3.2 in \cite{CP}. 
 
There are two types of $\Omega$, the one depending on  $P$ and the other not.  
The latter case is simple. $e^{ Q}$ acts on $\hOmega'$ as 
$e^{ Q}.\hOmega'=\Omega\otimes e^{ Q}\xi$. The finiteness and irreducibility of $\hV$  
imply the existence of a unique non-zero vector $\xi$ such that $e^Q\xi=C\xi$ with a 
complex number $C\not=0$.  Redefining $\frac{1}{C}e^Q$ as $e^Q$, we identify $\xi$ with $1$.  

For $\Omega$ depending on $P$, let us write the $P$ dependence explicitely as  
$\hOmega'(P)=\Omega(P)\otimes \xi$.  
$e^{ Q}$ acts on $\hOmega'(P)$ as $e^{ Q}.\hOmega'(P)=\Omega(P+ 1)\otimes e^{ Q}\xi$. 
The finiteness of $\hV$ implies that a finite number of vectors in 
$\{\Omega(P+n)\ (n\in \Z)\}$ are $\FF$-linearly independent. 
Setting $\widehat{\Omega}=\sum_{n\in \Z}\hOmega(P+n)\otimes \xi$, 
we have $e^{ Q}.\widehat{\Omega}=\sum_{n\in \Z}\hOmega(P+n)\otimes e^{ Q}\xi$. 
Then the same argument as the first 
case leads to $\xi=1$, and we obtain $\hOmega$ satisfying $3)$. 

In both cases, Theorem\ref{PBW} yields $\hV=U'_{q,p}.\hOmega$.
As for the action of $q^c$ on $\hV$, the statement follows from Corollary 3.2  
in \cite{CP}. \qed

\noi
{\it Remark.} An example of the vector $\hOmega'$ independent of $P$ is $v^l_0\otimes 1$ in 
Theorem \ref{repDri}, whereas the one depending on $P$ is $v^{(s)}$ in Theorem \ref{singvec}.

\begin{dfn}[Elliptic loop algebra]
The elliptic loop algebra $U_{q,p}(L(\slt))$ is the quotient of  
$U'_{q,p}(\slth)$ by the two sided ideal generated by $c$. 
\end{dfn}
Note 
 $U_{q,p}(L(\slt))\cong \FF[U_q(L(\slt))]\otimes_{\C}\C[\bH^*] $, where $\FF[U_q(L(\slt))]$ denotes the quantum loop algebra obtained as the quotient of $\FF[U_q(\slth)]$ by the two sided ideal generated by $c$\cite{CP}.  Note also that $U_{q,p}(L(\slt))$ 
is an $\h$-Hopf algebroid with the same $\mu_l, \mu_r, \Delta, \vep, S$ as 
$U'_{q,p}(\slth)$. 
Furthermore 
the $RLL$ relation for $U_{q,p}(L(\slt))$
 is given by \eqref{RLL3} with replacing $R^{+*}(u,P)$ with $R^{+}(u,P)$. 
It is identified  with 
the one for Felder's elliptic quantum group studied in \cite{EV1,TV,KNR}. 
Hence the corresponding $\hL^+(u)$ 
in  \eqref{Gauss} gives a realization of Felder's elliptic quantum group in terms of 
$U_{q,p}(L(\slt))$. 

Hereafter we 
consider dynamical representations of $U_{q,p}(L(\slt))$. 
\begin{dfn}[Pseudo-highest weight representation]\lb{hwrep}
A  dynamical representation $(\widehat{\pi}_V,\hV=V\otimes V_Q)$ of $U_{q,p}(L(\slt))$ is said to be pseudo-highest weight, if there exists a vector (pseudo-highest weight vector) $\hOmega\in \hV$ such that  $e^{ Q}.\hOmega=\hOmega$ and 
\be
&1)&x^+_n.\hOmega=0\quad (n\in \Z)\qquad\\
&2)&\psi_n.\hOmega=d^+_n\hOmega,\qquad \phi_{-n}.\hOmega=d^-_{-n}\hOmega\quad (n\in \Z_{\geq 0}),
\\
&3)&\hV=U_{q,p}(L(\slt)).\hOmega,
\en
with some complex numbers $d^\pm_{\pm n}$ satisfying $d^+_0d^-_0=1$.  
We call the set {\boldmath $d$}$=\{d^\pm_{\pm n}\}_{n\in \Z_{\geq0}}$ the pseudo-highest weight. 
\end{dfn}
We can state the equivalent conditions in terms of the matrix elements of $\hL^+(u)$. 
\begin{thm}\lb{hwrep2}
For a vector $\hOmega\in \hV$ satisfying $e^{ Q}.\hOmega=\hOmega$, 
the conditions $1)$ and $2)$ in Definition \ref{hwrep} are equivalent to the following. 
\be
&i)&\gamma(u).\hOmega=0\quad \forall u, \\
&ii)&q^h.\hOmega=q^\la\hOmega \quad \exists \la\in \C,\\
&& \al(u).\hOmega=A(u)\hOmega, 
\qquad \delta(u).\hOmega=D(u)\hOmega
\en
with some meromorphic functions $A(u)$ and $D(u)$ satisfying 
$D(u-1)^{-1}=A(u)$ and
\bea
&&A(u)=z^{\frac{\la}{2r}}\sum_{m\in \Z,n\in \Z_{\geq0}}
A_{m,n}z^mp^n \qquad A_{m,n}\in \C,\ z=q^{2u},\ p=q^{2r}.\lb{expandA}
\ena
\end{thm}
\noi
{\it Proof.} 
 We show that $i)$ and $ii)$ yield $1)$ and $2)$. 
Let us define $e_n\ (n\in \Z)$ by
\be
&&\phi_r(x^+(z))=\sum_{n\in \Z}e_n z^{-n}.
\en
From \eqref{Eplus}, we have \cite{JKOS2}  
\be
E^+(u)=e^{2Q}a^*[1]\sum_{n\in\Z}e_n\frac{1}{1-q^{2(P-1)}p^{n}}z^{-n-\frac{P-1}{r}}. 
\en
Here we used the following formula.  
\be
\frac{[u+s]}{[u][s]}=-\sum_{n\in \Z}\frac{1}{1-q^{-2s}p^n}.
\en
Then it follows from \eqref{Lelements} that $i)$ is equivalent to $ e_n.\hOmega=0$ 
for all $n\in \Z$. 

Furthermore from the definition of $x^+(z)$ and \eqref{uplus}, we have 
\be
&&e_n=\sum_{k\in\Z_{\geq 0}}p_k\left(\left\{\frac{a_{-l}q^{rl}}{[rl]_q}\right\}\right)x^+_{n+k}.
\en  
 Here $p_k(\{\al_l\})$ denotes the Schur polynomial defined by 
\be
&&\exp\left\{\sum_{n\in \Z_{>0}}\al_n z^n\right\}=\sum_{k\in \Z_{\geq 0}}p_k(\{\al_l\})z^{k}.
\en
$p_k(\{\al_l\})$ has the following expression.  
\be
&&p_k(\{\al_l\})=\sum_{m_1+2m_2+\cdots+km_k=k}\frac{\al_1^{m_1}
\cdots \al_k^{m_k}}{m_1!\cdots m_k!}.
\en 
Expanding $p_k\left(\left\{\frac{a_{-l}q^{rl}}{[rl]_q}\right\}\right)$ 
as a power series in $p=q^{2r}$, it follows that the condition 
$e_n.\hOmega=0$\ for all $n\in \Z$ is equivalent to $ x^+_n.\hOmega=0$\ for all $n\in \Z$. 

Similarly applying $ii)$,  we can evaluate $H^+(u).\hOmega$ as follows. 
\be
&&(q^{r}z)^{\frac{\la}{r}}u^+(z,p)\psi(z)u^-(z,p)
.\hOmega= A(u+1)A(u)\hOmega. 
\en
Here we used Proposition \ref{hpm} in the LHS, and \eqref{HKK} and \eqref{Lelements} in the RHS.
 Note that due to \eqref{expandA} fractional powers of $z$ in the both hand sides cancel out 
each other.  Expanding the both sides as a Laurent series in $z$ and 
a power series in $p$, one finds that $a_k\ (k\in \Z_{\not=0})$ are 
simultaneously diagonalized on $\hOmega$ and their eigenvalues are 
determined by the coefficients of the series in the right hand side.  \qed

\begin{dfn}[Verma module]
Let {\boldmath $d$}$=\{d^\pm_{\pm n}\}_{n\in \Z_{\geq0}}$
 be any sequence of complex numbers. 
The Verma module $M(\mbox{\boldmath $d$})$ is the quotient of  
$U_{q,p}(L(\slt))$ by the left ideal generated by $\{x^+_k\ (k\in \Z),\ \psi_n-d^+_n\cdot 1, \phi_{-n}-d^-_{-n}\cdot 1\ (n\in \Z_{\geq 0})
,\  e^{ Q}-1
\}$.
\end{dfn}
\begin{prop}\lb{verma}
The Verma module $M(\mbox{\boldmath $d$})$ is a pseudo-highest weight representation of 
 pseudo-highest weight {\boldmath $d$}. Every  pseudo-highest weight representation with  pseudo-highest weight {\boldmath $d$} is isomorphic to 
a quotient of $M(\mbox{\boldmath $d$})$. Moreover $M(\mbox{\boldmath $d$})$ has a unique 
maximal proper submodule $N(\mbox{\boldmath $d$})$, and up to isomorphism, $M(\mbox{\boldmath $d$})/N(\mbox{\boldmath $d$})$ is the unique irreducible pseudo-highest weight module of $U_{q,p}(L(\slt))$. 
 \end{prop}

\subsection{Elliptic Analogue of the Drinfeld Polynomials}\lb{classificationth}
We now consider a classification of finite-dimensional 
irreducible dynamical representations of $U_{q,p}(L(\slt))$. 
 We introduce a natural elliptic analogue of the Drinfeld polynomials. 
\begin{thm}\lb{eDriPoly}
The irreducible pseudo-highest weight dynamical 
 representation $(\widehat{\pi}_V,\hV)$ of $U_{q,p}(L(\slt))$ is 
 finite-dimensional if and only if  
there exists an  
entire and quasi-periodic function $P_V(u)$  such that 
\be
&&H^\pm(u).\hOmega=c_V\frac{P_V(u+1)}{P_V(u)}\hOmega,\\
&&P_V(u+r)=(-)^{{\rm deg}P}P_V(u),\\
&&P_V(u+r\tau)=(-)^{{\rm deg}P}e^{-\pi i\sum_{j=1}^{{\rm deg}P}(\frac{2(u-\al_j)}{r}+\tau)}P_V(u).
\en 
Here $\hOmega$ denotes the pseudo-highest weight vector in $\hV$, and $\tau=-\frac{2\pi i}{\log p}$. The symbol $c_V$ denotes a constant given by
\be
&&c_V=q^{\frac{r-1}{r}{\rm deg}P}\prod_{j=1}^{{\rm deg}P} a_j^{\frac{1}{r}},
\en
where ${\rm deg}P$ is a number of zeros of $P_V(u)$ in the fundamental parallelogram $(1,\tau)$ (= the degree of the Drinfeld polynomial $P(z)=\lim_{r\to \infty}P_V(u),\ z=q^{2u}$), and 
$a_j=q^{2\al_j}$ with $\al_j$ being a zero of $P_V(u)$ in the fundamental parallelogram.
The function $P_V(u)$ is unique up to a scalar multiple. 
\end{thm}
\noi
{\it Proof of the ``only if" part.}  From Theorem \ref{FDI}, $\hV$ has the 
pseudo-highest weight vector $\hOmega$. 
From Theorem 3.4 in \cite{CP}, there exists the Drinfeld 
polynomial $P(z)$ such that $P(0)=1$ and 
\be
&&\varphi(z).\hOmega=q^{{\rm deg}P}\frac{P(q^{-2}z^{-1})}{P(z^{-1})}\hOmega=\psi(z).\hOmega.
\en
Here the first and second equalities are in the sense of
 the power series in $z$ and $z^{-1}$, respectively.  
Then using Proposition \ref{hpm} and the formulae
\be
&&u^+(z,p)=\prod_{l=0}^\infty q^{h}\varphi(q^{{c}/{2}}q^{2r^*(l+1)}z),\quad 
u^-(z,p)=\prod_{l=0}^\infty q^{-h}\psi(q^{{c}/{2}}q^{-2r(l+1)}z),\lb{umpsi}
\en
we 
obtain 
\be
H^+(u).\hOmega&=&(q^rz)^{\frac{h}{r}}
\prod_{l=0}^\infty q^h\varphi(q^{2r(l+1)}z)\cdot \psi(z)\cdot
 \prod_{l=0}^\infty q^{-h}\psi(q^{-2r(l+1)}z)
\hOmega\\
&=&(q^rz)^{\frac{h}{r}}\prod_{l=1}^\infty 
q^{{\rm deg}P}\frac{P(q^{-2}q^{-2rl}z^{-1})}{P(q^{-2rl}z^{-1})}
 \prod_{l=0}^\infty q^{{\rm deg}P}\frac{P(q^{-2}q^{2rl}z^{-1})}{P(q^{2rl}z^{-1})}\hOmega.
\en
Supposing that the Drinfeld polynomial $P(z)$ is factorized as 
$P(z)=\prod_{j=1}^{{\rm deg}P}(1-a_jz)$, we have 
\be
 H^+(u).\hOmega&=&(q^rz)^{\frac{h}{r}}q^{{\rm deg}P}\prod_{j=1}^{{\rm deg}P}
\frac{\Theta_{q^{2r}}(a_j/q^2z)}{\Theta_{q^{2r}}(a_j/q^2z)}\hOmega\\
&=&q^{\frac{r-1}{r}{\rm deg}P}
\prod_{j=1}^{{\rm deg}P}a_j^{\frac{1}{r}}\frac{[u+1-\al_j]}{[u-\al_j]}\hOmega.
\en
This is the desired result with $P_V(u)=\prod_{j=1}^{{\rm deg}P}[u-\al_j]$. 
The quasi-periodicity of $P_V(u)$ follows from the one of the theta function $[u]$. 

The proof of the ``if" part is given in the next subsection. 
\qed

\noi
{\it Remark.} We can take $c_V=1$ by the gauge transformation
 given from (2.11) in \cite{JKOS2}.  An example is given in 
Corollary \ref{evP}.

\begin{prop}\lb{PVPW}
Let $\hV$ and $\widehat{W}$ be finite dimensional dynamical representations of 
$U_{q,p}(L(\slt))$ and assume that the 
tensor product $\hV\tot \widehat{W}$ is irreducible. Let 
$P_V(u), P_W(u)$ and $P_{V\tot W}(u)$ be the entire quasi-periodic 
function associated to $\hV,\hW$ and $\hV\tot \hW$ in Theorem \ref{eDriPoly}. Then
\be
&&P_{V\tot W}(u)=P_V(u)P_W(u).
\en
\end{prop}
\noi
{\it Proof.} The statement follows from the comultiplication formulae for the half currents in 
Proposition \ref{coprohalfc}. \qed

\subsection{Evaluation Representations} 
We consider an elliptic analogue of the evaluation representation of 
 $U_{q}(L(\slt))$\cite{Jimbo, CP}. This is an important example of the 
finite-dimensional irreducible dynamical representation of 
$U_{q,p}(L(\slt))$.  Some formulae presented here were essentially  
obtained in \cite{JKOS2}. Corollary \ref{evP} and Proposition \ref{LandR} are new. 

Let us consider the $l+1$-dimensional evaluation representation $(\pi_{l,w},V^{(l)}_w)$ of 
$\FF[U_q(L(\slt))]$. 
Here $V^{(l)}=\oplus_{m=0}^l\FF v^l_m,\ V^{(l)}_w=V^{(l)}\otimes \C[w,w^{-1}]$, and we define operators $h, S^\pm$ on $V^{(l)}$ by
\be
h v^{l}_m =(l-2m)v^l_m,\qquad S^{\pm}v^l_m=v^l_{m\mp1}, \qquad v^l_m=0\quad {\rm for } \ m<0,\ \ m>l.
\en
The action of the Drinfeld generators on $V^{(l)}_w$ is given as follows. 
\bea
&&\pi_{l,w}(a_n)=\frac{w^n}{n}\frac{1}{q-q^{-1}}((q^n+q^{-n})q^{nh}-(q^{(l+1)n}+q^{-(l+1)n})),\lb{evaluationrep}\\
&&\pi_{l,w}(x^{\pm}(z))=S^{\pm}\left[\frac{\pm h+l+2}{2}\right]_{q}\delta\left(q^{h\pm1}\frac{w}{z}\right).\nn
\ena

Applying this to Proposition \ref{dynamicalrep} and 
noting Definition \ref{halfcurrents}, we obtain the following theorem.  
\begin{thm}\lb{repDri}
Let $\hV^{(l)}(w)=V^{(l)}(w)\otimes \C1$ be the vector space, 
on which $e^{Q}$ acts as 
\be
&&e^{ Q}.(f(P)v\otimes 1)=f(P+ 1)v\otimes 1.
\en
The image of the half currents by 
the map $\widehat{\pi}_{l,w}={\pi}_{l,w}\otimes \id$ on $U_{q,p}(L(\slt))\cong \FF[U_q(L(\slt))]
\otimes_{\C} \C[\bH^*]$ is given, up to fractional powers of $z, w$ and $q$,   by 
\be
\widehat{\pi}_{l,w}(K^+(u))&=&-\frac{\varphi_l(u-v)}{[u-v-\frac{h-1}{2}]}e^{Q},\\
\widehat{\pi}_{l,w}(E^+(u))&=&-e^QS^+\frac{[u-v-\frac{h+1}{2}-P][\frac{l+h+2}{2}]}{[u-v-\frac{h+1}{2}][P]}e^Q,\\
\widehat{\pi}_{l,w}(F^+(u))&=&S^-\frac{[u-v+\frac{h-1}{2}+P][\frac{l-h+2}{2}]}{[u-v-\frac{h-1}{2}][P+h-1]},\\
\widehat{\pi}_{l,w}(H^{\pm}(u))&=&\frac{[u-v-\frac{l+1}{2}][u-v+\frac{l+1}{2}]}{
[u-v-\frac{h-1}{2}][u-v-\frac{h+1}{2}]}e^{2Q}.
\en
where $z=q^{2u}$, $w=q^{2v}$, and  
\be
\varphi_l(u)&=&-z^{-\frac{l}{2r}}\rho_{1l}^+(z,p)^{-1}[u+\frac{l+1}{2}],\\
\rho_{kl}^+(z,p)&=&q^{\frac{kl}{2}}\frac{\{pq^{k-l+2}z\}\{pq^{-k+l+2}z\}}{\{pq^{k+l+2}z\}\{pq^{-k-l+2}z\}}\frac{\{q^{k+l+2}/z\}\{q^{-k-l+2}/z\}}{\{q^{k-l+2}/z\}\{q^{-k+l+2}/z\}}.
\en
Furthermore $(\widehat{\pi}_{l,w}, \hV^{(l)}(w))$ is the $l+1$-dimensional 
irreducible dynamical representation of $U_{q,p}(L(\slt))$ with the  pseudo-highest weight 
vector $v_0^l\otimes 1$. 
\end{thm}

\noi
{\it Proof.} 
One can directly check that $\widehat{\pi}_{l,w}(K^+(u)), \widehat{\pi}_{l,w}(E^+(u))$ 
and $\widehat{\pi}_{l,w}(F^+(u))$ satisfy the relations 
in Theorem \ref{relHalfcurrents}. In the process, we use the formula
\bea
&&\varphi_l(u)\varphi_l(u-1)=[u-\frac{l+1}{2}][u+\frac{l+1}{2}].\lb{phiphi}
\ena
\qed

From Definition \ref{Loperator}, we obtain the image of the matrix elements of $\hL^+(u)$ 
as follows.
\begin{thm}\lb{repL}
\be
\widehat{\pi}_{l,w}(\al(u))&=&-\frac{[u-v+\frac{h+1}{2}][P-\frac{l-h}{2}][P+\frac{l+h+2}{2}]}{\varphi_l(u-v)[P][P+h+1]}e^Q,\\
\widehat{\pi}_{l,w}(\beta(u))&=&-S^-\frac{[u-v+\frac{h-1}{2}+P][\frac{l-h+2}{2}]}{\varphi_l(u-v)[P+h-1]}e^{-Q},\\
\widehat{\pi}_{l,w}(\gamma(u))&=&S^+\frac{[u-v-\frac{h+1}{2}-P][\frac{l+h+2}{2}]}{\varphi_l(u-v)[P]}e^Q,\\
\widehat{\pi}_{l,w}(\delta(u))&=&-\frac{[u-v-\frac{h-1}{2}]}{\varphi_l(u-v)}e^{-Q}.
\en
\end{thm}

\begin{cor}\lb{evP}
The elliptic analogue of the Drinfeld polynomial associated to $\hV^{(l)}(q^{2v})$ is given by
\be
&&P_{l,v}(u)=[u-v-\frac{l-1}{2}][u-v-\frac{l-1}{2}+1]\cdots 
[u-v+\frac{l-1}{2}].
\en
\end{cor}
\noi
{\it Proof.}
Noting $\delta(u)=K^+(u)^{-1}$, from \eqref{HKK}, \eqref{kplush} and Theorem \ref{repL}, we obtain 
\bea
&&\widehat{\pi}_{l,w}(H^\pm(u))(v^l_0\otimes 1)=\frac{[u-v+\frac{l+1}{2}]}{[u-v-\frac{l-1}{2}]}v^l_0\otimes 1.\lb{evHpm}
\ena
Then the entireness and the quasi-periodicity of $P_{l,v}(u)$ yield the desired result. 
\qed

Note that the zeros of $P_{l,v}(u)$ coincides with those of the 
Drinfeld polynomial corresponding to the evaluation representation $V^{(l)}(q^{2v})$ of 
 $U'_q(L(\slt))$ modulo $\Z r+\Z r\tau$. Note also that we have no $c_V$ factor in \eqref{evHpm} due to the remark below Theorem \ref{eDriPoly}.

{\it Proof of the ``if" part of Theorem \ref{eDriPoly}.} 
Let $P_V(u)$ be any entire quasi-periodic function satisfying the conditions 
in the Theorem \ref{eDriPoly}, and let its zeros in the 
fundamental parallelogram $(1,\tau)$ be 
$\al_1,\cdots, \al_r$. From Theorem \ref{hwrep2}, $P_V(u)$ determines the set of eigenvalues {\boldmath $d$} of $\psi_k$ and $\phi_{-k}\ (k\in \Z_{\geq 0})$ uniquely.  Consider the representation $\hV=\hV^{(1)}(q^{2\al_1})\tot \cdots \tot \hV^{(1)}(q^{2\al_1})$. Let 
${v}^1_0(i)=v^1_0\otimes 1$ denote the pseudo-highest weight vector in 
$\hV^{(1)}(q^{2\al_i})$ and set $\hOmega={v}^1_0(1) \tot \cdots \tot {v}^1_0(r) $. 
Then 
up to a scalar multiple, $\hOmega$ is a unique pseudo-highest weight vector such that $q^h.\hOmega=q^r\hOmega$. Let us consider the submodule  $\hV'=U_{q,p}(L(\slt)).\hOmega$ of $\hV$. 
 $\hV'$ has a unique maximal submodule $\hV''$. Then the quotient module $\hV'/\hV''$ is irreducible, and from Corollary \ref{evP} and Theorem \ref{PVPW}, $\hV'/\hV''$ has the entire quasi-periodic function given by
\be
&\widetilde{P_V}(u)=\prod_{j=1}^r[u-\al_j]. 
\en
$\widetilde{P_V}(u)$ has the same quasi-periodicity and zeros as $P_V(u)$. 
Hence $\widetilde{P_V}(u)$ coincides with $P_V(u)$ up to a scalar multiple. \qed

The following Proposition indicates a consistency of our construction of $\widehat{\pi}_{l,w}$ 
 and the fusion construction of the 
dynamical $R$ matrices (=face type Boltzmann weights).  
\begin{prop}\lb{LandR}
Let us define the matrix elements of $\widehat{\pi}_{l,w}(\hL^+_{\vep_1\vep_2}(u))$ by
\be
\widehat{\pi}_{l,w}(\hL^+_{\vep_1\vep_2}(u))v^l_m&=&
\sum_{m'=0}^l (\hL^+_{\vep_1\vep_2}(u))_{\mu_{m'}\mu_m}v^l_{m'},
\en
where $\mu_m=l-2m$. Then we have
\be
(\hL^+_{\vep_1\vep_2}(u))_{\mu_{m'}\mu_{m}}=
R^+_{1l}(u-v,P)_{\vep_1 \mu_{m'} }^{\vep_2 \mu_m}.
\en
Here $R^+_{1l}(u-v,P)$ is the $R$ matrix from (C.17) in  \cite{JKOS2}. 
The  case $l=1$, $R^+_{11}(u-v,P)$ coincides with the image 
$(\pi_{1,z}\otimes \pi_{1,w})$ of the universal $R$ matrix $\cR^+(\la)$\cite{JKOS} given in \eqref{Rmat}.
The case $l>1$, $R^+_{1l}(u-v,P)$ coincides with the $R$ matrix 
obtained by fusing $R^+_{11}(u-v,P)$ $l$-times. In particular the matrix element $R^+_{1l}(u-v,P)_{\vep \mu}^{\vep'\mu'}$ is gauge equivalent to the fusion face weight $W_{l1}(P+\vep',P+\vep'+\mu',P+\mu,P|u-v)$ from (4) in \cite{DJMO}. 
\end{prop}

\subsection{Tensor Product Representations}

In this subsection, we investigate a submodule structure of the tensor product space 
 $\hV^{(l_1)}(q^{2a})\tot \hV^{(l_2)}(q^{2b})$ and derive an elliptic analogue of the 
Clebsch-Gordan coefficients. 
We abbreviate the pseudo-highest weight vectors ${v}^{l_1}_0
\otimes 1$ and 
${v}^{l_2}_0\otimes 1$ of $\hV^{(l_1)}(q^{2a})$ and $\hV^{(l_2)}(q^{2b})$ as  
${v}^{l_1}_0$ and 
${v}^{l_2}_0$, respectively. 

\begin{thm}\lb{singvec}
There exists a  vector $v^{(s)}\in
\hV^{(l_1)}(q^{2a})\tot \hV^{(l_2)}(q^{2b})$ satisfying the conditions $1)\sim 3)$ in the below, 
if and only if  $b-a=\frac{l_1+l_2-2s}{2}+1\ (s=0,1,\cdots,{\rm min}\{l_1,l_2\})$.  
\be
          &&1)\ q^h.v^{(s)}=q^{l_1+l_2-2s}v^{(s)}\quad (h\in\h),\quad\\
          &&2)\ \Delta(\gamma(u)).v^{(s)}=0\quad \forall u,\\
          &&3)\ \Delta(\al(u)).v^{(s)}= A(u)v^{(s)}, \qquad 
            \Delta(\delta(u)).v^{(s)}= D(u)v^{(s)}\qquad \forall u, 
\en
where
\be 
A(u)&=&\frac{[u-a-\frac{l_1+1}{2}][u-a+\frac{l_1+1}{2}]}{\varphi_{l_1}(u-a)\varphi_{l_2}(u-b)},\\
D(u)&=&\frac{[u-a-\frac{l_1-1}{2}+s][u-a-\frac{l_1-1}{2}-l_2+s-1]}{\varphi_{l_1}(u-a)\varphi_{l_2}(u-b)}.
\en
Explicitely, the vector $v^{(s)}$ is given by
\bea
v^{(s)}&=&\sum_{m_1=0}^{s}C^s_{m_1}(P)v^{l_1}_{m_1}\tot v^{l_2}_{s-m_1},\nn\\
C_{m_1}^s(P)&=&C_{0}^s \frac{[P-l_2+s-m_1]_{s-m_1}[l_2-s+1]_{m_1}}
{[P+1]_{s-m_1}[-l_1]_{m_1}}.\lb{L:coeff}
\ena
\end{thm} 
\noi
{\it Proof.}  We solve the conditions $1) \sim 3)$. 
The first condition yields
\bea
v^{(s)}&=&\sum_{m_1=0}^{s}C^s_{m_1}(u,P)v^{l_1}_{m_1}\tot v^{l_2}_{s-m_1}\lb{vscvv}
\ena
with unknown coefficients $C^s_{m_1}(u,P)$. 

By using $2)$ and \eqref{shifts}, we obtain   
\be
\Delta(\gamma(u))v^{(s)}&=&\sum_{m_1}\left\{
C^s_{m_1}(u,P+1)\gamma(u)v^{l_1}_{m_1}\tot \al(u)v^{l_2}_{s-m_1}
+C^s_{m_1}(u,P-1)\delta(u)v^{l_1}_{m_1}\tot \gamma(u)v^{l_2}_{s-m_1}
\right\}\\
&=&0.
\en
Apply Theorem \ref{repL} and move all the coefficients in the second tensor space to the first one 
by using the following formula obtained from \eqref{VtotW} 
\bea
v\tot f(u,P)v^{l_2}_{m_2}
&=&v\tot  f(u,P+h-(l_2-2m_2))v^{l_2}_{m_2}\nn\\
&=& f(u,P-(l_2-2m_2))v\tot v^{l_2}_{m_2}.\lb{totf}
\ena
Here one should note $f^*(u,P)=f(u,P)$, i.e. $c=0$, in the evaluation representations. 
We thus obtain the following recursion relation.
\bea
C^s_{m_1}(u,P)&=&-C^s_{m_1-1}(u,P-2)\frac{[u-a-\frac{l_1+1}{2}+m_1][u-b
+\frac{l_2-2s+1}{2}+m_1+1-P]}{[u-a-\frac{l_1+1}{2}+m_1+1-P][u-b
+\frac{l_2-2s+1}{2}+m_1]}\nn\\
&&\quad \times\frac{[l_2-s+m_1][P][P-1]}{[l_1+1-m_1][P-l_2+s-1-m_1][P+s-m_1]}.\lb{recursionc}
\ena
Then one finds that all the $u$ dependent factors cancel out each other, 
if and only if $b-a=\frac{l_1+l_2-2s}{2}+1\ (s=0,1,\cdots,{\rm min}\{l_1,l_2\})$. 
We hence obtain 
\bea
C^s_{m_1}(P)&=&-C^s_{m_1-1}(P-2)
\frac{[l_2-s+m_1][P][P-1]}{[l_1+1-m_1][P-l_2+s-1-m_1][P+s-m_1]}.\lb{recursionc2}
\ena
Here we rewrote $C^s_{m_1}(u,P)$ as $C^s_{m_1}(P)$. 
Solving this, we obtain 
\bea
C^s_{m_1}(P)&=&C^s_{0}(P-2m_1)\frac{[l_2-s+1]_{m_1}[P-2m_1+1]_{2m_1}}
{[-l_1]_{m_1}[P-l_2+s-2m_1]_{m_1}[P+s-2m_1+1]_{m_1}}.\lb{coeffv}
\ena
Here $[u]_m$ denotes the elliptic shifted factorial 
\be
[u]_{m}&=&[u][u+1]\cdots [u+m-1].
\en

Finally the third condition yields, for $\delta(u)$, 
\be
&&\hspace*{-5mm}\Delta(\delta(u))v^{(s)}\\
&&=\sum_{m_1=0}^sC_{m_1}^s(P+1)\gamma(u)v^{l_1}_{m_1}\tot \beta(u)v^{l_2}_{s-m_1}
+\sum_{m_1=0}^sC_{m_1}^s(P-1)\delta(u)v^{l_1}_{m_1}\tot \delta(u)v^{l_2}_{s-m_1}\\
&&=\sum_{m_1=0}^s\frac{C_{m_1}^s(P-1)}{\varphi_{l_1}(u-a)\varphi_{l_2}(u-b)
[P-l_2+s-m_1-1][P+s-m_1]}\\
&&\hspace*{-10mm}\times \left(
[P-l_2+s-m_1-1][P+s-m_1][u-a-\frac{l_1-2m_1-1}{2}][u-b-\frac{l_2-2(s-m_1)-1}{2}]\right.\\
&&\hspace*{-10mm}\left. +[s-m_1][l_2-s+m_1+1][u-a-\frac{l_1-2m_1-1}{2}-P][u-b-\frac{l_2-2(s-m_1)-1}{2}+P]\right)v^{l_1}_{m_1}\tot v^{l_2}_{s-m_1}. 
\en
In the second equality, we used \eqref{recursionc2}. 
By using the theta function identity
\be
[u+x][u-x][v+y][v-y]-[u+y][u-y][v+x][v-x]=[x-y][x+y][u+v][u-v], 
\en
we obtain
\be
\Delta(\delta(u))v^{(s)}&=&D(u) \sum_{m_1=0}^s C^s_{m_1}(P-1)\frac{[P][P-l_2+2s-2m_1-1]}{
[P-l_2+s-m_1-1][P+s-m_1]}v^{l_1}_{m_1}\tot v^{l_2}_{s-m_1}
\en
with $D(u)$ appearing in the statement of the theorem. 
Then the necessary and sufficient condition that $\Delta(\delta(u))$ is 
diagonalized on $v^{(s)}$ with the eigen value $D(u)$ is 
\be
C_{m_1}^s(P)=C^s_{m_1}(P-1)\frac{[P][P-l_2+2s-2m_1-1]}{
[P-l_2+s-m_1-1][P+s-m_1]}.
\en
Solving this, we obtain
\bea
C_{m_1}^s(P)=C_{m_1}^s\frac{[P-l_2+s-m_1]_{s-m_1}}{[P+1]_{s-m_1}}.\lb{coeffb}
\ena
Here $C_{m_1}^s$ is a coefficient independent of $P$, which 
 can be determined by substituting \eqref{coeffb} into \eqref{recursionc2}. We hence obtain 
$C_{m_1}^s(P)$ in the form in \eqref{L:coeff}. 
We can check that for $\Delta(\al(u))$, the similar argument leads to the same result 
\eqref{L:coeff}.
\qed

\noi
{\it Remark.} A similar statement was obtained in \cite{FV}.

By applying $\beta(u)$ on $v^{(s)}$ repeatedly, we can compute the other weight vectors as follows. 

\begin{thm}\lb{thmbetavs}
Setting $l=l_1+l_2-2s$, we have for  $0\leq m\leq l$ 
\bea
&&\Delta(\beta(u)\beta(u+1)\cdots\beta(u+m-1))v^{(s)}\nn\\
&&=\frac{[P]}{\prod_{i=1}^m\varphi_{l_1}(u-a+i-1)\varphi_{l_2}(u-a-
\frac{l}{2}+i-1)}\nn\\
&&\times\hspace*{-10mm}
\sum_{k={\rm max}(0,s+m-l_2)}^{{\rm min}(l_1,s+m)}
\hspace*{-6mm}(-)^kC_0^s\frac{[P+m-2k-l_2+s]_s}{[P+m-2k+1]_s}
[u-a+\frac{l_1+1}{2}]_{m-k}
[u-a-l+\frac{l_1-1}{2}+m-k+P]_{m-k}
\nn\\
&&\qquad\times
[-u+a-\frac{l_1-1}{2}-m+k-P]_{k}[-u+a+l-\frac{l_1-1}{2}-m+1]_k \nn\\
&&\qquad\times\frac{[-m]_k
[P-k]_{m-k}[P+l_1-k+1]_{m-k}[s+1]_{m-k}}
{[P-k+1]_{m-k}[P]_{m-k}[P+l_1-2k+1]_m}\nn\\
&&\qquad\times 
{}_{12}V_{11}\left(P+m-2k;-k,-s,P-k,l_2-s+1,-u+a-\frac{l_1-1}{2},\right.\nn\\
&&\qquad\qquad \left. u-a-l+
\frac{l_1-1}{2}+2m-2k+P,P+m-2k+l_1+1\right)
\ v^{l_1}_k \tot\ v^{l_2}_{m+s-k}.\lb{betavs}
\ena
Moreover, we have   for $l< m$
\be 
\Delta(\beta(u)\beta(u+1)\cdots\beta(u+m-1))v^{(s)}=0.
\en
\end{thm}
Here ${}_{12}V_{11}$ denotes the very-well-poised balanced elliptic hypergeometric series  defined by\cite{FT,SZ, S} 
\be
{}_{s+1}V_s(u_0;u_1,\cdots,u_{s-4})=\sum_{j=0}^\infty\frac{[u_0+2j]}{[u_0]}
\prod_{i=0}^{s-4}\frac{[u_i]_j}{[u_0+1-u_i]_j}
\en
with the balancing condition
\be
\sum_{i=1}^{s-4}u_i=\frac{s-7}{2}+\frac{s-5}{2}u_0.
\en
We give the proof in Appendix \ref{proofA}.

The vector $v^{(s)}$ obtained in Theorem \ref{singvec} depends on the dynamical parameter $P$. 
Let us write its $P$ dependence explicitely as 
$v^{(s)}(P)$.  It satisfies 
$e^{ Q}.v^{(s)}(P)=v^{(s)}(P+ 1)$.  Setting $\widehat{v}^{(s)}=\sum_{n\in \Z}v^{(s)}(P+n)$, 
we have $e^{ Q}.\widehat{v}^{(s)}=\widehat{v}^{(s)}$. Hence $\widehat{v}^{(s)}$ is a 
 pseudo-highest weight vector. 
Note that $\widehat{v}^{(s)}$ is a vector in the $\FF$-linear space spaned by 
the vectors $v^{l_1}_{m_1}\tot v^{l_2}_{s-m}\ (0\leq m_1\leq s)$.  
Let us consider the pseudo-highest weight 
$U_{q,p}(L(\slt))$-module $W^{(s)}$ generated by $\widehat{v}^{(s)}$.

\begin{thm}\lb{submodule}
If $b-a=\frac{l_1+l_2-2s}{2}+1\quad 0< s\leq {\rm min}(l_1,l_2)$, 
the pseudo-highest weight $U_{q,p}(L(\slt))$-module  $W^{(s)}$ is a unique 
proper submodule of $\hV=\hV^{(l_1)}(q^{2a})\tot 
\hV^{(l_2)}(q^{2b})$. Moreover we have 
\bea
W^{(s)}\cong \hV^{(l_1-s)}(q^{2(a-\frac{s}{2})})\tot 
\hV^{(l_2-s)}(q^{2(b+\frac{s}{2})}),
\ena
\bea
\hV/W^{(s)}\cong \hV^{(s-1)}(q^{2(a+\frac{l_1-s+1}{2})})\tot \hV^{(l_1+l_2-s+1)}(q^{2(b-\frac{l_1-s+1}{2})}).
\ena
\end{thm}

\noi
{\it Proof.} 
From Theorem \ref{eDriPoly}, it is enough to show that the entire quasi-periodic functions 
associated with the representations in the both sides coincide with each other under the condition 
$b-a=\frac{l_1+l_2-2s}{2}+1$ for $0<s\leq {\rm min}(l_1,l_2)$. 
In fact, one can evaluate the action of the operators $H^{\pm}(u)$ 
on the highest weight vectors $\widehat{v}^{(s)}$ and 
$v^{l_1-s}_0\tot v^{l_2-s}_0$ of $W^{(s)}$ and 
$\hV^{(l_1-s)}(q^{2(a-\frac{s}{2})})\tot 
\hV^{(l_2-s)}(q^{2(b+\frac{s}{2})})$, respectively, 
and find that the eigen values coincide with each other. They are
 given by 
 \be
&&\frac{[u-a-\frac{l_1+1}{2}][u-a+\frac{l_1+1}{2}]}{
[u-a-\frac{l_1-1}{2}+s][u-a-\frac{l_1+1}{2}-1+s]}. 
\en

Similarly,  
the eigen values of $H^{\pm}(u)$ 
on the highest weight vectors $v^{l_1}_0\tot v^{l_2}_0$ and 
$v^{s-1}_0\tot v^{l_1+l_2-s+1}_0$ of $\hV/W^{(s)}$ and 
$\hV^{(s-1)}(q^{2(a+\frac{l_1-s+1}{2})})\tot \hV^{(l_1+l_2-s+1)}(q^{2(b-\frac{l_1-s+1}{2})})$, respectively, coincide and 
are given by 
 \be
&&\frac{[u-a+\frac{l_1+1}{2}][u-b+\frac{l_2+1}{2}]}{
[u-a-\frac{l_1-1}{2}][u-b-\frac{l_2-1}{2}]}. 
\en
\qed

\noi
{\it Remark.} A similar statement was presented in \cite{FV} without proof.

\section{Discussions}
In this section, we consider the limits  trigonometric $r\to \infty$, non-affine  
$u\to \infty$ and  non-dynamical  $P\to \infty$ of the results and  make some remarks on
their algebraic structures and relations.  
  
In the limits, the elliptic dynamical 
$R$ matrix degenerates  as follows. 
\be
\mmatrix{
R^+(u,P)&\stackrel{r\to\infty}{\longrightarrow}& R^+_{trig.}(u,P)&
 \stackrel{u\to\infty}{\longrightarrow}& R^+(P)& \stackrel{P\to\infty}{\longrightarrow}& R^+,\cr
        &                                      &                 &
{\searrow \atop P\to\infty}&  & {\nearrow \atop u\to\infty} & \\[-4mm]
        &                                      &                  & &R^+(u)& &\cr}  
\en
where setting $x=q^{2P}$ we have   
\be
R^+_{trig.}(u,P)&=&\rho^+_{trig.}(z)
\mat{1&0&0&0\cr
   0&\frac{(1-q^2x)(1-q^{-2}x)}{(1-x)^2}b(z)&\frac{1-xz}{1-x}c(z)&0\cr
0&\frac{1-xz^{-1}}{1-x}zc(z)&b(z)&0\cr
0&0&0&1\cr},\\
R^+(P)&=&q^{1/2}\mat{1&0&0&0\cr
           0&\frac{q(1-q^2x)(1-q^{-2}x)}{(1-x)^2}&
\frac{1-q^2}{1-x}&0\cr
0&-\frac{x(1-q^2)}{1-x}&q&0\cr
0&0&0&1\cr},\quad 
R^+(u)=\rho^+_{trig.}(z)
\mat{1&0&0&0\cr
   0&b(z)&c(z)&0\cr
0&zc(z)&b(z)&0\cr
0&0&0&1\cr},
\en
\be
 R^+&=&q^{1/2}\mat{1&0&0&0\cr
    0&q&1-q^2&0\cr
0&0&q&0\cr
0&0&0&1\cr},
\en
\be
\rho^+_{trig.}(z)&=&q^{1/2}\frac{(z^{-1};q^4)_\infty
(q^4z^{-1};q^4)_\infty}{(q^2z^{-1};q^4)^2_\infty},\\
b(z)&=&\frac{q(1-z)}{1-q^2z},\qquad c(z)\ =\ \frac{1-q^2}{1-q^2z}.
\en

Correspondingly naive degeneration limits of the $RLL$ relation \eqref{RLL3} 
imply the following diagram of quantum algebras. 
\bea
\mmatrix{U_{q,p}(\slth)&\stackrel{r\to\infty}{\longrightarrow}& U_{q,x}(\slth)&
 \stackrel{u\to\infty}{\longrightarrow}& U_{q,x}(\slt)& \stackrel{P\to\infty}{\longrightarrow}&
 U_{q}(\slt).\cr
        &                                      &                 &
{\searrow \atop P\to\infty}&  & {\nearrow \atop u\to\infty} & \\[-4mm]
        &                                      &                  & &U_q(\slth)& &\cr}  
\lb{degQG}
\ena   
Here $U_{q,x}(\slth)$ is the dynamical quantum affine algebra suggested 
in \cite{AAFRR}. However neither its generators nor the $L$ operators 
has yet been given explicitely. Let us make some speculations on it.  
The $U_{q,x}(\slth)$ should be a semi-direct product $\C$-algebra isomorphic to 
$\FF[U_q(\slth)]\otimes_{\C} \C[\bH^*]$ 
and characterized by the $RLL$ relation of the type \eqref{RLL3} 
associated with  $R^+_{trig.}(u,P)$. 
In fact,  the $\hL^+(u)$ operator as well as 
the half currents in Definition \ref{halfcurrents} do survive 
in the trigonometric limit. The $\h$-Hopf algebroid structure 
of $U_{q,p}(\slth)$ also survives in the limit. 
We hence expect that  $U_{q,x}(\slth)$ is an $H$-Hopf algebroid. 
We will discuss this subject in elsewhere.  

$U_{q,x}(\slt)$ denotes  the dynamical quantum algebra introduced by Babelon
\cite{Babelon}, and  $U_q(\slth)$,  $U_q(\slt)$ are the standard   
  quantum  affine and non-affine algebras by Drinfeld-Jimbo, respectively. 
The quasi-Hopf algebra structure of $U_{q,x}(\slt)$ was studied in \cite{BBB}, whereas 
the generalized FRST formulation and the $\h$-Hopf algebroid structure  
were discussed in \cite{EV1,EV2,KR}. The FRST formulation and Hopf algebra structure of 
$U_q(\slth)$ and  $U_q(\slt)$ were given in \cite{RS} and \cite{FRT}, respectively. 

Concerning the FRST formulations, we should remark that 
in the elliptic algebra $U_{q,p}(\slth)$ as well as in $\Bqla(\slth)$ 
the $L$ operators $\hL^+(u)$ and $\hL^-(u)$ are not independent \cite{JKOS}.
This is also true for its trigonometric and non-affine limits. 
However this is not true 
after the non-dynamical limit $P\to \infty$, so that  we need 
two $L$ operators $L^+$ and $L^-$ for $U_q(\slth)$ and $U_q(\slt)$.

It is also interesting to see the limits of  
${}_{12}V_{11}$ obtained in \eqref{betavs}. Corresponding to the upper 
series in  \eqref{degQG}, we find the following.
\bea
&&{}_{12}V_{11}\left(P+m-2k;-s,-k,P-k,l_2-s+1,-u+a-\frac{l_1-1}{2},\right.\nn\\
&&\qquad\qquad \left. u-a-l+
\frac{l_1-1}{2}+2m-2k+P,P+m-2k+l_1+1\right) \nn\\
&\stackrel{r\to\infty}{\longrightarrow}&
{}_{10}W_{9}\left(q^{2(P+m-2k)};q^{-2s},q^{-2k},q^{2(P-k)},q^{2(l_2-s+1)},q^{2(-u+a-\frac{l_1-1}{2})},\right.\nn\\
&&\qquad\qquad \left. q^{2(u-a-l+
\frac{l_1-1}{2}+2m-2k+P)},q^{2(P+m-2k+l_1+1)};q^2,q^2\right)\nn\\
& \stackrel{u\to\infty}{\longrightarrow}&
{}_{8}W_{7}\left(q^{2(P+m-2k)};q^{-2s},q^{-2k},q^{2(P-k)},q^{2(l_2-s+1)},q^{2(P+m-2k+l_1+1)};q^2,q^{-2(l-m)}\right)\nn\\
&&\qquad 
=\frac{(q^{2(P+m-2k+1)},q^{2(m+1)},q^{2(P+m-k-l_2+s)},q^{2(k-l_1)};q^2)_s}{
(q^{2(P+m-k+1)},q^{2(m-k+1)},q^{2(P+m-2k-l_2+s)},q^{-2l_1};q^2)_s}q^{-2sk}\nn\\
&&\qquad\qquad  \times {}_4\phi_3
\left(\mmatrix{q^{-2s},q^{-2k},q^{-2(P+m-k+s)},q^{2(l-m+1)}\cr
q^{-2(s+m)},q^{-2(P+m-k+l_1-l-1)},q^{2(l_1+1-s-k)}\cr};q^2,q^2\right)\nn\\
& \stackrel{P\to\infty}{\longrightarrow} &
\frac{(q^{2(m+1)};q^2)_s}{
(q^{2(m-k+1)};q^2)_s}
 {}_3\phi_2
\left(\mmatrix{q^{-2s},q^{-2k},q^{-2(s+l+1)}\cr
q^{-2(s+m)},q^{-2l_1}\cr};q^2,q^2\right),\lb{limits12V11}
\ena
where 
\be
&&(a_1,a_2,\cdots,a_m;q^2)_s=\prod_{i=1}^m(a_i;q^2)_s,\\
&&(a;q^2)_s=(1-a)(1-aq^2)\cdots(1-aq^{2(s-1)}).
\en
Here we followed the notations in \cite{GR}. 
After the second limit, we used the transformation 
formula from (2.17) in \cite{KR} 
whereas  after the third limit, the formula from (3.2.2) in \cite{GR}.  
As shown by Rosengren \cite{RoEle},  ${}_{10}W_{9}$ 
in the above gives a system of biorthogonal 
functions identical to the one obtained by Wilson \cite{Wilson}. 
The ${}_4\phi_3$ part is identified with the $q$-Racah  polynomial, 
and the ${}_3\phi_2$ part with 
the $q$-Hahn polynomial. 

A representation theoretical derivation of ${}_3\phi_2$ or the $q$-Clebsch-Gordan coefficients 
was done on the basis of $U_q(\slt)$ in \cite{KirRes,KoeKoo,Vak}. 
  The case of ${}_8W_7$ or the $q$-Racah polynomials, or Askey-Wilson polynomial, has an 
interesting history. It was first done on the basis of the quantum group 
$SU_{q}(2)$ with considering the so-called twisted primitive element 
in \cite{Koorn} and \cite{NM}. 
Later an alternative derivation was carried out on the basis of the co-representations of 
$U_{q,x}(\slt)$ in \cite{KR}. 
The relation between these two derivations can be found in \cite{Ro3} and \cite{Stokman}.  
As for the case ${}_{12}V_{11}$,    
a representation theoretical derivation was first done 
in  \cite{KNR} on the basis of co-representations of 
Felder's elliptic quantum group.  In this paper we have given an alternative derivation 
on the basis of representations of $U_{q,p}(\slth)$. 
Comparing \eqref{degQG} and \eqref{limits12V11}, we conjecture that 
Wilson's biorthogonal functions ${}_{10}W_{9}$ can be 
derived similarly on the basis of $U_{q,x}(\slth)$.  

Moreover it is instructive to note that reading the diagram \eqref{limits12V11} in inverse 
direction the dynamical parameter $P$ modifies the $q$-$3j$-symbols  
(${}_3\phi_2$ or the Clebsch-Gordan coefficients) into the $q$-$6j$-symbols (${}_4\phi_3$ or 
$q$-Racah polynomials), 
the affinization parameter $u$ modifies the orthogonal polynomials ($q$-Racah, $q$-Hahn
 polynomials ) into the biorthogonal 
functions (Wilson's biorthogonal function).
As a result ${}_{12}V_{11}$ is an elliptic analogue of the 
$q$-$6j$-symbol and is biorthogonal. We think that 
this observation should become a guiding principle 
in choosing a suitable type of quantum groups, such as dynamical or non-dynamical, 
affine or non-affine, in a derivation of elliptic analogues of the $q$-special functions. 

It is also interesting to note that 
the above degeneration diagram of ${}_{12}V_{11}$ 
coincides with the one of the hypergeometric type special solutions of the discrete 
Painlev\'{e} equations \cite{KMNOYell} corresponding to the following degeneration of  
 affine Weyl group symmetries \cite{Sakai}. 
\be
&&E^{(1)}_8 \to E^{(1)}_7\to E^{(1)}_6\to D^{(1)}_5\to\cdots.
\en 
It should be interesting if one could find a direct 
connection between this diagram or the discrete Painlev\'{e} equations themselves 
and  the 
quantum groups in \eqref{degQG}. 

\subsection*{Acknowledgments}

The author would like to thank Michio Jimbo, Anatol Kirillov, Atsushi Nakayashiki,  Masatoshi Noumi, 
Hjalmar Rosengren and Tadashi Shima 
for stimulating discussions and valuable suggestions. 
He is also grateful to Tetsuo Deguchi, Jonas Hartwig, Masahiko Ito, Masaki Kashiwara, 
Christian Korf, Barry McCoy,  Tetsuji Miwa, Tomoki Nakanishi, Masato Okado, 
Vitaly Tarasov and Yasuhiko Yamada for their interests and 
useful conversations. He also thank Hjalmar Rosengren for his kind 
hospitality during a stay in Charmers University of Technology and 
G\"{o}teborg University. This work is supported by the Grant-in-Aid for 
Scientific Research (C)19540033, JSPS
Japan.

\appendix

\setcounter{equation}{0}
\begin{appendix}

\section{The $RLL$ relation (2.17) at $c=0$}
We write down the $RLL$ relation \eqref{RLL3} in terms of the matrix elements of 
$\hL^+(u)$ in the case $c=0$. 
\bea
&&[\al(u_1),\al(u_2)]=0,\qquad [\delta(u_1),\delta(u_2)]=0,\\
&&[\beta(u_1),\beta(u_2)]=0,\qquad [\gamma(u_1),\gamma(u_2)]=0,\\
&&\al(u_1)\beta(u_2)=\bar{c}(u,P)\al(u_2)\beta(u_1)+b(u,P)\beta(u_2)\al(u_1),\lb{commab}\\
&&\beta(u_1)\al(u_2)=\bar{b}(u,P)\al(u_2)\beta(u_1)+c(u,P)\beta(u_2)\al(u_1),\\
&&\gamma(u_1)\delta(u_2)=\bar{c}(u,P)\gamma(u_2)\delta(u_1)+b(u,P)\delta(u_2)\gamma(u_1),\\
&&\delta(u_1)\gamma(u_2)=\bar{b}(u,P)\gamma(u_2)\delta(u_1)+c(u,P)\delta(u_2)\gamma(u_1),
\\
&&{c}(u,P+h)\gamma(u_1)\al(u_2)+b(u,P+h)\al(u_1)\gamma(u_2)=\gamma(u_2)\al(u_1),\\
&&\bar{b}(u,P+h)\gamma(u_1)\al(u_2)+\bar{c}(u,P+h)\al(u_1)\gamma(u_2)=\al(u_2)\gamma(u_1),\\
&&{c}(u,P+h)\delta(u_1)\beta(u_2)+b(u,P+h)\beta(u_1)\delta(u_2)=\delta(u_2)\beta(u_1),\\
&&\bar{b}(u,P+h)\delta(u_1)\beta(u_2)+\bar{c}(u,P+h)\beta(u_1)\delta(u_2)=\beta(u_2)\delta(u_1),
\lb{commdb}\\
&&{c}(u,P+h)\gamma(u_1)\beta(u_2)+b(u,P+h)\al(u_1)\delta(u_2)\nn\\
&&\qquad=\bar{c}(u,P)\gamma(u_2)\beta(u_1)+b(u,P)\delta(u_2)\al(u_1),\\
&&\bar{b}(u,P+h)\gamma(u_1)\beta(u_2)+\bar{c}(u,P+h)\al(u_1)\delta(u_2)\nn\\
&&\qquad=b(u,P)\beta(u_2)\gamma(u_1)+\bar{c}(u,P)\al(u_2)\delta(u_1),\\
&&{b}(u,P+h)\beta(u_1)\gamma(u_2)+c(u,P+h)\delta(u_1)\al(u_2)\nn\\
&&\qquad=\bar{b}(u,P)\gamma(u_2)\beta(u_1)+c(u,P)\delta(u_2)\al(u_1),\\
&&\bar{c}(u,P+h)\beta(u_1)\gamma(u_2)+\bar{b}(u,P+h)\delta(u_1)\al(u_2)\nn\\
&&\qquad={c}(u,P)\beta(u_2)\gamma(u_1)+\bar{b}(u,P)\al(u_2)\delta(u_1).
\ena

\section{Proof of Theorem 4.18} \lb{proofA}
In order to prove the theorem, we need the following four Lemmas.
\begin{lem}\lb{albetas}
For $c=0$,
\be
\al(u)\beta(v_1)\cdots\beta(v_l)
&=&\frac{[P+1][P-l][u-v_l]}{[P][P-l+1][u-v_1+1]}\beta(v_1)\cdots\beta(v_l)\al(u)\nn\\
&&+\sum_{k=1}^l\frac{[P+1][P-k+1-u+v_k][1]}{[P][u-v_1+1][P-k+2]}\beta(v_1)\cdots \al(v_k)\cdots \beta(v_l)\beta(u).
\en
\end{lem}
\noi
{\it Proof.} Use \eqref{commab} repeatedly. 
\qed

\begin{lem}\lb{ellbi}
\be
&&\Delta(\beta(u)\beta(u+1)\cdots\beta(u+m-1))\\
&&\qquad=\sum_{j=0}^m D^m_j(P)\al(u+m-1)\cdots\al(u+m-j)\beta(u+m-j-1)\cdots\beta(u)\\
&&\qquad\qquad\qquad \tot\ \delta(u)\cdots\delta(u+m-j-1)
\beta(u+m-j)\cdots\beta(u+m-1),
\en
where
\bea
&&D^m_j(P)=\frac{[1]_m}{[1]_j[1]_{m-j}}\frac{[P][P-m+2j]}{[P+j][P-m+j]}\qquad  m\in \Z_{\geq 0}.\lb{ellbinomial}
\ena
\end{lem}

\noi
{\it Proof.} We prove the statement by induction on $m$. The case $m=1$ is just the 
comultiplication formula for $\beta(u)$. Assume that the statement is true for $m$. 
Then
\be
&&\Delta(\beta(u)\beta(u+1)\cdots\beta(u+m-1))\Delta(\beta(u+m))\\
&&=\sum_{j=0}^m D^m_{j}(P)\al(u+m-1)\cdots \al(u+m-j)\beta(u+m-j-1)\cdots\beta(u)\al(u+m)
\\
&&\qquad \tot \delta(u)\cdots\delta(u+m-j-1)\beta(u+m-j)\cdots\beta(u+m-1)\beta(u+m)\\
&&+\sum_{j=0}^m D^m_{j}(P)\al(u+m-1)\cdots \al(u+m-j)\beta(u+m-j-1)\cdots\beta(u)
\beta(u+m)\\
&&\qquad \tot \delta(u)\cdots\delta(u+m-j-1)\beta(u+m-j)\cdots\beta(u+m-1)\delta(u+m)\\
&&=\al(u+m-1)\cdots \al(u)\al(u+m)+\beta(u+m-1)\cdots\beta(u)\beta(u+m)\\
&&+\sum_{j=1}^m \left\{{\ \atop\ } D^m_{j-1}(P)
\al(u+m-1)\cdots \al(u+m-j+1)\beta(u+m-j)\cdots\beta(u)\al(u+m)\right.\\
&&\left.+\frac{[P+j-m]}{[P+2j-m]}D^m_{j}(P)\al(u+m-1)\cdots \al(u+m-j)\beta(u+m-j-1)
\cdots
\beta(u)\beta(u+m)\right\}\\
&&\tot \delta(u)\cdots\delta(u+m-j-1)\delta(u+m-j)\beta(u+m-j+1)\cdots\delta(u+m).
\en
To obtain the second equality we used the property of $\tot$ in \eqref{fstotf} with putting $c=0$ 
and the following relation obtained from \eqref{commdb} 
with putting $u_1=v$ and $u_2=v+1$  
\be
&&\delta(v)\beta(v+1)=\frac{[P+h+1]}{[P+h]}\beta(v)\delta(v+1). 
\en
Therefore we need to show
\bea
&&D^{m+1}_j(P-j+1)\al(u+m)\beta(u+m-j)\cdots\beta(u)\nn\\
&&=D^m_{j-1}(P-j+1)\beta(u+m-j)\cdots\beta(u)\al(u+m)\nn\\
&&+\frac{[P-m+1]}{[P+j-m+1]}D^m_{j}(P-j+1)\al(u+m-j)\beta(u+m-j-1)\cdots\beta(u)\beta(u+m)\lb{commab0}
\ena
for $j=1,2,\cdots,m$.

Specializing $l\to m-j, u\to u+m-j, v_k\to u+k-1\ (1\leq k\leq m-j)$ in Lemma \ref{albetas}, 
we have
\be
&&\al(u+m-j)\beta(u+m-j-1)\cdots\beta(u)\nn\\
&&=\sum_{k=1}^{m-j+1}\frac{[P+1][P-m+j][1]}{[P][m-j+1][P-k+2]}\beta(u)\cdots\al(u+k-1)\cdots\beta(u+m-j).
\en
Substituting this into the second term in the RHS of \eqref{commab0}, we obtain\bea
&&D^{m+1}_j(P-j+1)\al(u+m)\beta(u+m-j)\cdots\beta(u)\nn\\
&&=D^m_{j-1}(P-j+1)\beta(u+m-j)\cdots\beta(u)\al(u+m)\nn\\
&&\qquad+D^m_{j}(P-j+1)
\sum_{k=1}^{m-j+1}\frac{[P-m+1][P+1][P-m+j][1]}{[P+j-m+1][P][m-j+1][P-k+2]}\nn\\
&&\qquad\qquad \times\beta(u)\cdots\al(u+k-1)\cdots\beta(u+m-j)\beta(u+m)\lb{commab1}.
\ena

Similarly, specializing $l\to m-j+1, u\to u+m, v_k\to u+k-1\ (1\leq k\leq m-j+1)$ 
 in Lemma \ref{albetas}~, we obtain
\bea
&&\al(u+m)\beta(u+m-j)\cdots\beta(u)\nn\\
&&=\frac{[P+1][P-m+j-1][j]}{[P][P-m+j][m+1]}\beta(u)\cdots\beta(u+m-j)\al(u+m)\nn\\
&&+\sum_{k=1}^{m-j+1}\frac{[P+1][P-m][1]}{[P][m+1][P-k+2]}
\beta(u)\cdots\al(u+k-1)\cdots\beta(u+m-j)\beta(u+m).\lb{commab2}
\ena

We compare this with \eqref{commab1}. 
These two relations coincide with each other if and only if 
\be
&&
\frac{D^m_{j-1}(P-j+1)}{D^{m+1}_{j}(P-j+1)}=\frac{[P+1][P-m+j-1][j]}{[P][P-m+j][m+1]},\\
&&
\frac{D^m_{j}(P-j+1)}{D^{m+1}_{j}(P-j+1)}=\frac{[P-m][P-m+j+1][m-j+1]}{[P-m+j][P-m+1][m+1]}.
\en
Therefore we obtain
\be
\frac{D^m_{j}(P)}{D^{m}_{j-1}(P)}=\frac{[P-m+j-1][P-m+2j][m-j+1][P+j-1]}
{[P-m+j][P-m+2(j-1)][P+j][j]}.
\en
Solving this with the initial condition $D^m_0(P)=1$, we obtain \eqref{ellbinomial}. \qed

\begin{lem}\lb{ab}
For $v^{l_1}_{m_1}\in \hV^{(l_1)}(q^{2a})$, we have 
\be
&&\hspace*{-1cm}\al(u+m-1)\cdots\al(u+m-j)\beta(u+m-j-1)\cdots\beta(u)v^{l_1}_{m_1}\\
&&\hspace*{-1cm}=(-)^{k+m_1+m}\frac{[u-a+\frac{l_1+1}{2}]_{m-k}[P-k]_{m-k}
[P+l_1-k+1]_{m-k}[-u+a-m-\frac{l_1-1}{2}-P+k]_{k}[1]_k}
{\prod_{i=1}^m\varphi_l(u-a+i-1)\ [P]_{m-k}[P+l_1-2k+1]_{m}}\\
&&\qquad\qquad\qquad \times 
\frac{[-u+a-\frac{l_1+1}{2}+1]_{m_1}[P-2k+m]_{m_1}[P+l_1+m-2k+1]_{m_1}}
{[P+m-k]_{m_1}[u-a+m+\frac{l_1-1}{2}+P-2k+1]_{m_1}[1]_{m_1}}\ v^{l_1}_{k}.
\en
Here we set $k=m_1+m-j$.
\end{lem}

\noi
{\it Proof.} Applying Theorem\ref{repL}, we evaluate the LHS as  
\be
{\rm LHS}&=&(-)^m\prod_{i=1}^j\frac{[u-a+m-i+\frac{h+1}{2}][P+i-1-\frac{l_1-h}{2}][P+i-1+\frac{l_1+h+2}{2}]}{\varphi_l(u-a+m-i)[P+i-1][P+h+i]}\\
&&\qquad\qquad\qquad \times 
\prod_{i=1}^{m-j}\frac{[u-a+m+\frac{h-1}{2}+P-i+1][\frac{l_1-h+2}{2}-i]}
{\varphi_l(u-a+m-j-i)[P+h+j+i]}\ v^{l_1}_{m_1+m-j}\\
&=&(-)^{k+m_1}\frac{[-u+a-m-\frac{l_1+1}{2}+1+k]_{m_1+m-k}[P-k]_{m_1+m-k}
[P+l_1-k+1]_{m_1+m-k}}
{\prod_{i=1}^m\varphi_l(u-a+i-1)\ [P]_{m_1+m-k}}\\
&&\qquad\qquad\qquad \times 
\frac{[-u+a-m-\frac{l_1-1}{2}+k-P]_{k-m_1}[-k]_{k-m_1}}{[P+l_1-2k+1]_{m}}\ v^{l_1}_{k}.
\en
Then using the following formulae, we obtain the desired result. 
For $b, k, m_1 ,k-m_1\in \Z_{\geq 0}$, 
\be
&&[a]_{m_1+b}=[a]_{b}[a+b]_{m_1},\qquad 
\\&&
[a-b]_{m_1+b}=(-)^b[-a+1]_{b}[a]_{m_1},    \qquad \\
&&[a+k]_{k-m_1}=(-)^{m_1}\frac{[a+k]_k}{[-a-2k+1]_{m_1}}.\hspace*{7cm}\qed
\en

\begin{lem}\lb{db}
For $v^{l_2}_{s-m_1}\in \hV^{(l_2)}(q^{2b})$, we have 
\be
&&\hspace*{-1cm}v^{l_1}_k\ \tot\ \delta(u)\cdots\delta(u+m-j-1)\beta(u+m-j)
\cdots\beta(u+m-1)\ v^{l_2}_{s-m_1}\\
&&\hspace*{-1cm}=(-)^{m+k}\frac{[-u+b+\frac{l_2-1}{2}-m-s+1]_{k}
[u-b-\frac{l_2+1}{2}+m+s+1-k+P]_{m-k}[s+1]_{m-k}}
{\prod_{i=1}^m\varphi_{l_2}(u-b+i-1)[P-k+1]_{m-k}}\\
&&\hspace*{-1cm}
\qquad\times \frac{[u-b-\frac{l_2+1}{2}+2m+s+1-2k+P]_{m_1}[P-k+1]_{m_1}[-s]_{m_1}}
{[-u+b+\frac{l_2-1}{2}-m-s+1]_{m_1}[P+m-2k+1]_{2m_1}}
v^{l_1}_k\ \tot\  v^{l_2}_{m+s-k}.
\en
Here we set $k=m_1+m-j$.
\end{lem}

\noi
{\it Proof.}
Applying Theorem \ref{repL}, we evaluate the LHS as
\be
&&\hspace*{-1cm}v^{l_1}_k\ \tot\ 
(-)^m\prod_{i=1}^{m-j}\frac{[u-b+i-1-\frac{h-1}{2}]}{\varphi_{l_2}(u-b+i-1)}
\prod_{i=1}^{j}\frac{[u-b+\frac{h-1}{2}+P+i][\frac{l_2-h+2}{2}-i]}{
\varphi_{l_2}(u-b+m-j+i-1)[P+h-m+j+i]}v^{l_2}_{s-m_1+j}\\
&&\hspace*{-1cm}=v^{l_1}_k\ \tot\ (-)^{m_1+k}
\frac{[u-b-\frac{l_2-1}{2}+m+s-k]_{k-m_1}
[u-b-\frac{l_2+1}{2}+m+s+P+h+1-k]_{m_1+m-k}
}{\prod_{i=1}^m\varphi_{l_2}(u-b+i-1)[P+h+m_1-k+1]_{m_1+m-k}}\\
&&\qquad\qquad \times [-m-s+k]_{m_1+m-k}\
 v^{l_2}_{m+s-k}.
\en
Then the statement follows from \eqref{totf} and 
\be
&&[a-k]_{k-m_1}=(-)^{k+m_1}\frac{[-a+1]_k}{[-a+1]_{m_1}},\\
&&[a+m_1-k]_{m_1+m-k}=\frac{[a-k]_{m-k}[a+m-2k]_{2m_1}}{[a-k]_{m_1}}\qquad k,m,m_1\in 
\Z_{\geq 0},\ m\geq k\geq m_1.\qed
\en

{\it Proof of the first statement in Theorem \ref{thmbetavs}.}  
By using Lemma \ref{ellbi} and \eqref{vscvv}, we obtain  
\be
&&\hspace*{-1cm}{\rm LHS}\\
&&\hspace*{-1cm}=\sum_{m_1=0}^s\sum_{j=0}^m
D_j^m(P)C_{m_1}^{s}(P-m+2j)\al(u+m-1)\cdots\al(u+m-j)\beta(u+m-j-1)\cdots\beta(u)
v^{l_1}_{m_1}\ \\
&&\hspace*{-1cm}\quad \tot\ \delta(u)\cdots\delta(u+m-j-1)\beta(u+m-j)
\cdots\beta(u+m-1)\ v^{l_2}_{s-m_1}.
\en
Then using \eqref{L:coeff} and Lemma \ref{ab}, \ref{db}, and change the summation variable from $j$ to $k$ by $k=m_1+m-j$,  we obtain the desired result.
In the process, ${}_{12}V_{11}$ is identified with the part associated with 
the summation  with respect to $m_1$ over 
${\rm max}(0,k-m)\leq m_1\leq {\rm min}(k,s)$. 
There we also manipulate  \eqref{L:coeff} by the formula
\be
\frac{[P-l_2+s-m_1]_{s-m_1}}{[P+1]_{s-m_1}}=\frac{[P-l_2+s-2m_1]_s[P-2m_1+1]_{2m_1}}{[P-2m_1+1]_s[P-l_2+s-2m_1]_{m_1}[P+s-2m_1+1]_{m_1}}.
\en

{\it Proof of the second statement.}  Let us set $m=l_1+l_2-2s+1$. 
Then  $l_1-s+1\leq k\leq l_1$. We show that 
${}_{12}V_{11}$ part in \eqref{betavs} vanishes for $k=l_1-s+n\ (n=1,2,\cdots, s)$. 
In fact, substituting  $m=l_1+l_2-2s+1$ and  $k=l_1-s+n$, we find that the 
${}_{12}V_{11}$ part is reduced to 
\bea 
&&\hspace*{-1.5cm}\sum_{m_1=0}^s\frac{[P-l_1+l_2+1-2n+2m_1]}{[P-l_1+l_2+1-2n]}\nn\\
&&\hspace*{-1.5cm}\times\frac{[P-l_1+l_2+1-2n]_{m_1}[-l_1+s-n]_{m_1}[-s]_{m_1}[P-l_1+s-n]_{m_1}[l_2+1-s]_{m_1}[P+l_2+2-2n]_{m_1}}{[1]_{m_1}[P+l_2-s+2-n]_{{m_1}}[P-l_1+l_2+2+s-2n]_{m_1}[l_2-s+2-n]_{m_1}[P-l_1+s+1-2n]_{m_1}[-l_1]_{m_1}}\nn\\
&&\hspace*{-1.5cm}=\frac{[1-n]_s[-P+l_1-l_2-2s-1+2n]_s[l_1+l_2-2s+2]_s[-P-s+n]_s}{[l_1-s+1]_s
[-P-l_2-1+n]_s[l_2-s+2-n]_s[-P+l_1-2s+2n]_s}. \lb{reduced12V11}
\ena
In the last line we used the Jackson-Frenkel-Turaev summation formula\cite{FT}
\be
&&{}_{10}V_{9}(\beta-\gamma-s;-s,\alpha-\gamma,-\alpha-\gamma+1-s,\beta+\delta,\beta-\delta)
=\frac{[\gamma-\beta,\gamma+\beta,\al+\delta,\al-\delta]_s}
{[\al-\beta,\al+\beta,\gamma+\delta,\gamma-\delta]_s}.
\en
\eqref{reduced12V11} vanishes for $n=1,2,\cdots, s$. \qed

\end{appendix}

\end{document}